\documentclass{amsart}

\usepackage{graphicx}
\usepackage[all]{xy}
\usepackage{color}

\newtheorem{thrm}{Theorem}

\newtheorem{defn}[thrm]{Definition}
\newtheorem{rem}[thrm]{Remark}
\newtheorem{prop}[thrm]{Proposition}

\def \Dj{\mbox{\raise0.3ex\hbox{-}\kern-0.4em D}}

\newcommand{\Crit}{\operatorname{Crit}}
\newcommand{\Sym}{\operatorname{Sym}}
\newcommand{\Sp}{\operatorname{Sp}}
\newcommand{\ind}{\operatorname{ind}}
\newcommand{\Ind}{\operatorname{Ind}}
\newcommand{\im}{\operatorname{im}}

\begin{document}

\title[PSS isomorphism and conormal spectral invariants]{Piunikhin--Salamon--Schwarz isomorphisms and spectral invariants for conormal bundle}

\author{Jovana \Dj ureti\'c}

\address{Matemati\v{c}ki fakultet, Studentski trg 16, 11000
Belgrade, Serbia}

\email{jovanadj@matf.bg.ac.rs}
\thanks{This work is partially supported by Serbian Ministry of Education, Science and Technological Development project \#174034}

\begin{abstract}
We give a construction of Piunikhin--Salamon--Schwarz isomorphism between the Morse homology and the Floer homology generated by Hamiltonian orbits starting at the zero section and ending at the conormal bundle. We also prove that this isomorphism is natural in the sense that it commutes with the isomorphisms between the Morse homology for different choices of the Morse function and the Floer homology for different choices of the Hamiltonian. We define a product on the Floer homology and prove triangle inequality for conormal spectral invariants with respect to this product.
\end{abstract}
\maketitle


MS Classification: Primary 53D40, Secondary 53D12, 57R58, 57R17


\section{\textbf{Introduction and main results}}

Let $M$ be a compact smooth manifold. The cotangent bundle $T^*{M}$ of $M$ carries a natural symplectic structure $\omega=d\lambda$, where $\lambda$ is the Liouville form. Let
$$
\nu^*N=\{\alpha\in T_p^*M\,:\,p\in N,\,\alpha|_{T_pN}=0\}\subset T^*M,
$$
 be a conormal bundle of a closed submanifold $N\subseteq M$.  Let $H$ be a time-dependent smooth compactly supported Hamiltonian on $T^*M$ such that the intersection $\nu^*N\cap\phi_H^1(o_M)$ is transverse. Here, $\phi_H^t:T^*M\to T^*M$ denotes Hamiltonian flow of Hamiltonian vector field $X_H$. Floer chain groups $CF_*(o_M,\nu^*N:H)$ are $\mathbb{Z}_{2}-$vector spaces generated by the finite set $\nu^*N\cap\phi_H^1(o_M)$ (see~\cite{P} for more details). Floer homology $HF_*(o_M,\nu^*N:H)$ is defined as the homology group of $(CF_*(o_M,\nu^*N:H),\partial_F)$ where $\partial_F$ is a boundary operator$$\partial_F(x)=\sum_{y\in\nu^*N\cap\phi_H^1(o_M)}n(x,y;H)y,$$ and $n(x,y;H)$ is the (mod 2) number of solutions of a system
\begin{equation}\label{eq:1}
\left\{
\begin{array}{ll}
\frac{\partial u}{\partial s}+J(\frac{\partial u}{\partial t}-X_H(u))= 0,\\
u(s, 0)\in o_M, \; u(s, 1)\in \nu^*N, \\
u(-\infty,t)=\phi_H^t((\phi_H^1)^{-1})(x), \;
u(+\infty,t)=\phi_H^t((\phi_H^1)^{-1})(y), \\
x,y \in \nu^*N\cap\phi_H^1(o_M). \\
\end{array}
\right.
\end{equation}
This homology was introduced by Floer in~\cite{F1}, developed by Oh in~\cite{O1} and Fukaya, Oh, Ohta and Ono in most general case (see~\cite{FOOO}). For a convenience, these groups will be denoted by $HF_*(H)$. Although it is well known that these groups do not depend on $H$, we will keep $H$ in the notation, since in many practical applications it is useful to keep track on the Hamiltonian used in their definition.
For two regular pairs of parameters $(H^\alpha,J^\alpha)$ and $(H^\beta,J^\beta)$ the isomorphism between corresponding Floer homology groups
$$
S^{\alpha\beta}:HF_*(H^\alpha)\to HF_*(H^\beta)
$$ is induced by a chain homomorphism
$$
\sigma^{\alpha\beta}:CF_*(H^\alpha)\to CF_*(H^\beta),\,\,\,\sigma^{\alpha\beta}(x^\alpha)=\sum_{x^\beta}n(x^\alpha,x^\beta;H^{\alpha\beta})x^\beta,
$$
that counts the number $n(x^\alpha,x^\beta;H^{\alpha\beta})$ of solutions of a system
\begin{equation}\label{eq:2}
\left\{
\begin{array}{ll}
\frac{\partial u}{\partial s}+J^{\alpha\beta}(\frac{\partial u}{\partial t}-X_{H^{\alpha\beta}}(u))= 0,\\
u(s, 0)\in o_M, \; u(s, 1)\in \nu^*N ,\\
u(-\infty,t)=\phi_{H^\alpha}^t((\phi_{H^\alpha}^1)^{-1})(x^\alpha), \;
u(+\infty,t)=\phi_{H^\beta}^t((\phi_{H^\beta}^1)^{-1})(x^\beta), \\
x^\alpha \in \nu^*N\cap\phi_{H^\alpha}^1(o_M),\,x^\beta \in \nu^*N\cap\phi_{H^\beta}^1(o_M). \\
\end{array}
\right.
\end{equation}
Here $H^{\alpha\beta}_s$ and $J^{\alpha\beta}_s$ are $s$-dependent families such that for some $R>0$

\[H^{\alpha\beta}_s = \left\{
\begin{array}{l l}
  H^\alpha, & \quad \mbox{$s\le-R$}\\
  H^\beta, & \quad \mbox{$s\ge R$},\\ \end{array} \right. \]
\[J^{\alpha\beta}_s = \left\{
\begin{array}{l l}
  J^\alpha, & \quad \mbox{$s\le-R$}\\
  J^\beta, & \quad \mbox{$s\ge R$}.\\ \end{array} \right. \]
\indent We define the action functional ${\mathcal A}_H$ on the space of paths
$$\Omega(o_M,\nu^*N)=\{\gamma:[0,1]\to T^*M\,|\,\gamma(0)\in o_M,\,\gamma(1)\in\nu^*N\}$$
by
$$
{\mathcal A}_H(\gamma)=-\int\gamma^*\lambda+\int_0^1H(\gamma(t),t)\,dt.
$$
Critical points of ${\mathcal A}_H$ are Hamiltonian paths with ends on the zero section and the conormal bundle, i.e. $CF_*(H)$. Now we can define filtered Floer homology. Denote by
$$
CF^\lambda_*(H)={\mathbb Z}_2\langle x\in CF_*(H)\,|\,{\mathcal A}_H(x)<\lambda\rangle.
$$
Since the action functional decreases along holomorphic strip (see \cite{O1} for details) differential $\partial_F$ preserves the filtration given by ${\mathcal A}_H$. Its restriction
$$
\partial_F^\lambda=\partial_F|_{CF^\lambda_*(H)}
$$
defines a boundary operator on the filtered complex $CF^\lambda_*(H)$. Filtered Floer homology is now defined as a homology of the filtered complex
$$
HF^\lambda_*(H)=H_*(CF^\lambda_*(H),\partial_F^\lambda).
$$
Note that filtered Floer homology depends on the Hamiltonian $H$.
\\
\indent Let us recall the definition of Morse homology. For a Morse function
$f:N\rightarrow{\mathbb R}$ Morse chain complex, $CM_*(N:f)$, is a ${\mathbb Z}_2$--vector space generated by the set of critical points of $f$. Morse homology
groups $HM_*(N:f)$ are the homology groups of $CM_*(N:f)$ with respect
to a boundary operator
$$
\partial_M:CM_*(N:f)\rightarrow CM_*(N:f), \quad
\partial_M (p)=\sum_{q\in \Crit (f)}n(p,q;f)q,
$$
where $n(p,q;f)$ is the number of gradient trajectories that satisfy
\begin{equation}\label{eq:3}
\left\{
\begin{array}{ll}
\frac{d\gamma}{ds}=-\nabla f(\gamma), \\
\gamma(-\infty)=p, \; \gamma(+\infty)=q.
\end{array}
\right.
\end{equation}
In a way analogous to $S^{\alpha\beta}$, we can define isomorphism between Morse homologies of two different Morse functions $f^{\alpha}$ and $f^{\beta}$
$$T^{\alpha\beta}:HM_*(f^\alpha)\rightarrow HM_*(f^\beta).$$
It is generated by a chain homomorphism
$$
\tau^{\alpha\beta}:CM_*(f^\alpha)\to CM_*(f^\beta),\,\,\,\tau^{\alpha\beta}(p^\alpha)=\sum_{p^\beta}n(p^\alpha,p^\beta;f^{\alpha\beta})p^\beta,
$$
that counts the number $n(p^\alpha,p^\beta;f^{\alpha\beta})$ of solutions of a system
\begin{equation}\label{eq:parmorse}
\left\{
\begin{array}{ll}
\frac{d\gamma}{ds}=-\nabla f^{\alpha\beta}(\gamma),\\
\gamma(-\infty)=p^{\alpha}, \; \gamma(+\infty)=p^{\beta},
\end{array}
\right.
\end{equation}
(see~\cite{Sc1} for details). We use brief notation $HM_*(f)$ instead of $HM_*(N:f)$. Morse homology groups $HM_*(f)$ are isomorphic to singular homology groups $H_*(N;{\mathbb Z}_2)$~\cite{M,Sa1,Sc1} (we will sometimes identify Morse and singular homologies).\\

\indent Our first theorem gives isomorphisms between Morse homology $HM_*(N:f)$ and Floer homology $HF_*(o_M,\nu^*N:H)$. These isomorphisms are essentially different from ones defined in~\cite{P}.
\begin{thrm}\label{main_thm}
There exist isomorphisms
$$\Phi:HF_k(o_M,\nu^*N:H)\rightarrow HM_k(N:f),
$$
$$
\Psi:HM_k(N:f)\rightarrow HF_k(o_M,\nu^*N:H),
$$
that are inverse to each other:
\begin{equation}\label{eq:isomor}
\Phi\circ\Psi={\mathbb {Id}}|_{HM} \quad \text{and} \quad \Psi\circ\Phi=\mathbb{Id}|_{HF}.
\end{equation}
\end{thrm}
\indent In order to obtain isomorphisms on homology level we consider homomorphisms on chain complexes defined by counting the intersection number of the space of gradient trajectories of function $f$ and the space of perturbed holomorphic discs with boundary on the zero section $o_M$ and the conormal bundle $\nu^*N$.
\vskip5mm
\begin{center}\includegraphics[width=8cm,height=4cm]{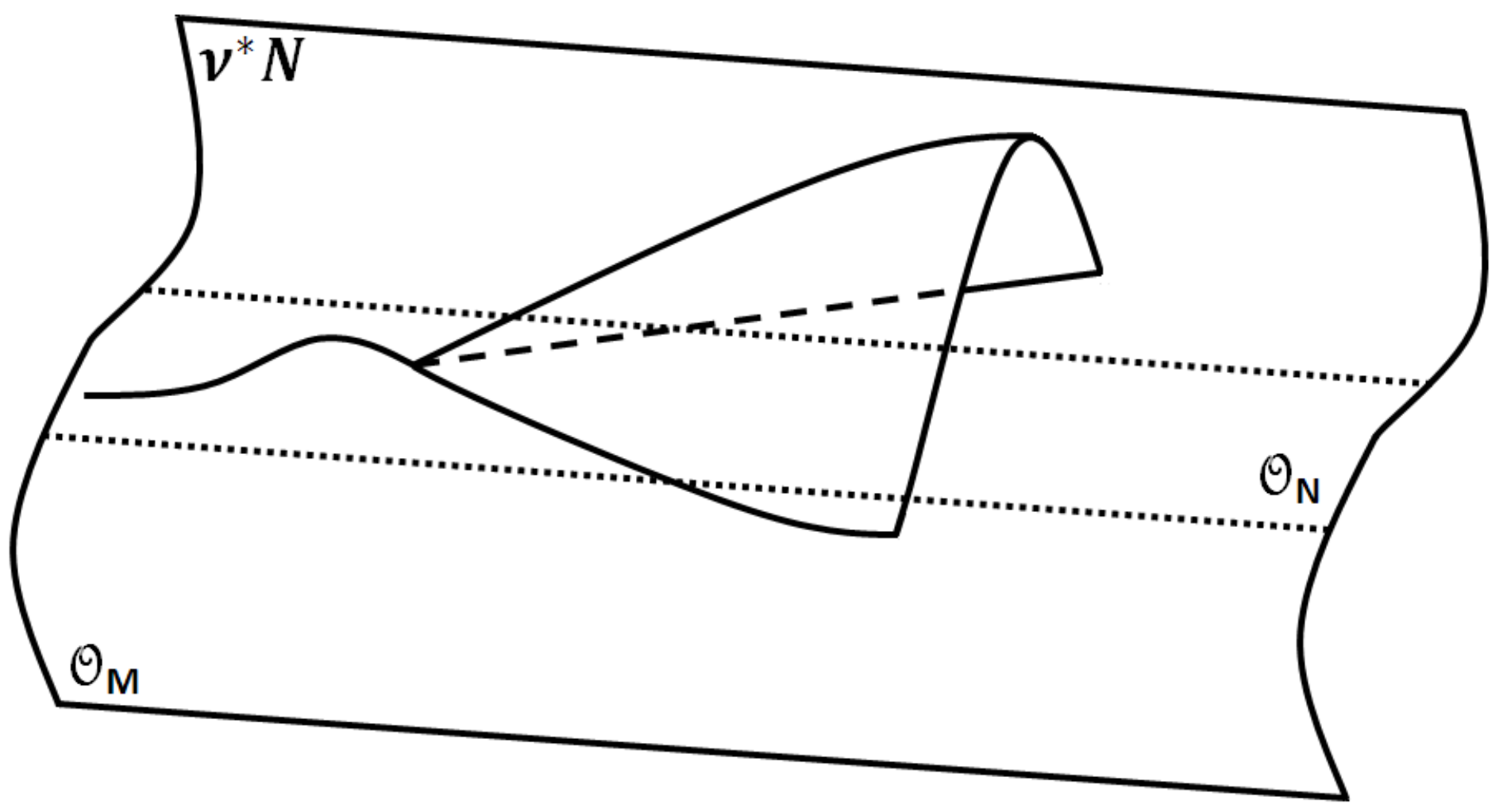}\\
Fig. 1. \textit{Intersection of gradient trajectory and perturbed holomorphic disc}
\end{center}
\vskip5mm

The main problem we need to overcome is that we have singular Lagrangian boundary conditions on holomorphic diks since an intersection $o_{M}|_{N}=o_M\cap\nu^*N$ is not transverse.\\
\indent Motivation for this isomorphism was a paper by Piunikhin, Salamon and Schwarz,\\~\cite{PSS} and a paper by Kati\'{c} and Milinkovi\'{c},~\cite{KM}. In~\cite{PSS} they considered Floer homology for periodic orbits. Kati\'{c} and Milinkovi\'{c} gave a construction of Piunikhin--Salamon--Schwarz isomorphisms in Lagrangian intersections Floer homology for a cotangent bundle. They worked with Floer homology generated by Hamiltonian orbits that start and end on zero section $o_M$. We obtain that isomorphism as special case for $N=M$. Albers constructed PSS--type homomorphism in more general symplectic manifold (which is not necessary an isomorphism, see~\cite{Alb}).\\
\indent In~\cite{P} Po\'{z}niak constructed a different type of isomorphism between Morse homology $HM_*(N:f)$ and Floer homology $HF_*(o_M,\nu^*N:H_f)$. Namely, he used Hamiltonian $H_f$ that is an extension of a Morse function $f$. We don't have that kind of restriction, our Hamiltonian $H$ doesn't have to be an extension of a Morse function $f$.\\
\indent Another advantage of using our isomorphism is its naturality. Using Po\'{z}niak's type isomorphism it is not obvious whether this diagram
\begin{equation}\label{moj_diagram}
\xymatrix{ HF_*(H^{\alpha})
 \ar[r]^-{S^{\alpha\beta}}
&HF_*(H^{\beta})\\
HM_*(f^{\alpha})
 \ar[r]^-{T^{\alpha\beta}}
 \ar[u]
&HM_*(f^{\beta})
 \ar[u]
,}\end{equation}
commutes because different type of equations are used in definitions of $S^{\alpha\beta}$ and $T^{\alpha\beta}$. If we use our, PSS--type, isomorphisms as vertical arrows we obtain commutativity of diagram~(\ref{moj_diagram}).
\begin{thrm}\label{diagram_commutes}
Diagram
\begin{equation}\label{konacan_dijagram}
\xymatrix{ HF_k(o_M,\nu^*N:H^{\alpha})
 \ar[r]^-{S^{\alpha\beta}}
&HF_k(o_M,\nu^*N:H^{\beta})\\
HM_k(N:f^{\alpha})
 \ar[r]^-{T^{\alpha\beta}}
 \ar[u]_{\Psi^\alpha}
&HM_k(N:f^{\beta})
 \ar[u]_{\Psi^\beta}
,}\end{equation}
commutes.
\end{thrm}
\indent Using the existence of PSS isomorphism we can define conormal spectral invariants and prove some of their properties. Denote by $$\imath^\lambda_*:HF_*^\lambda(H)\to HF_*(H)$$
the homomorphism induced by the inclusion map
$$\imath^\lambda:CF_*^\lambda(H)\to CF_*(H).$$
For $\alpha\in HM_*(N:f)$ define a conormal spectral invariant
$$l(\alpha;o_M,\nu^*N:H)=\inf\{\lambda\,|\,\Psi(\alpha)\in\im(\imath^\lambda_*)\}.$$
Oh defined Lagrangian spectral invariants in \cite{O1} using the idea of Viterbo's invariants for generating functions (see \cite{V}). It turns out that these two invariants are the same (under some normalizaton conditions), see \cite{M1,M2}.\\
\indent By counting a pair--of--pants with appropriate boundary conditions we prove that there exists a product on homology $HF_*(o_M,\nu^*N:H)$ (see figure 2).
\begin{thrm}\label{prod_thm}
Let $H_1,H_2,H_3\in C^{\infty}_c([0,1]\times T^*M)$ be three Hamiltonians with compact support. Then, there exists a product in homology
$$
\ast:HF_*(o_M,\nu^*N:H_1)\otimes HF_*(o_M,\nu^*N:H_2)\to HF_*(o_M,\nu^*N:H_3).
$$
The product $\ast$ induces the operation on $H_*(N)$ via the PSS isomorphism
$$
\alpha\cdot\beta=\Phi(\Psi(\alpha)\ast\Psi(\beta)),
$$
for $\alpha,\beta\in H_*(N)$.
\end{thrm}
\vskip2mm
\begin{center}\includegraphics[width=6cm,height=3cm]{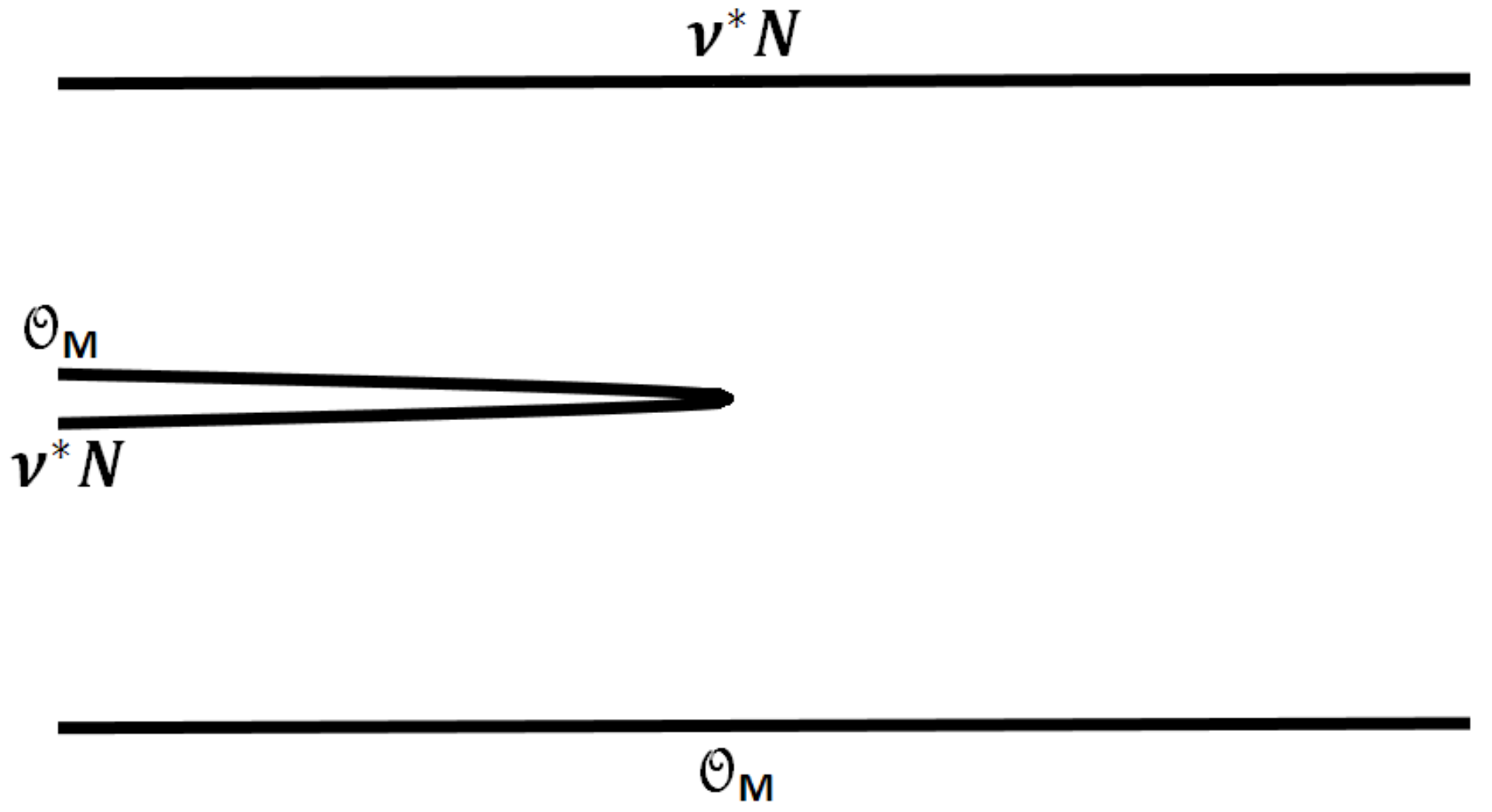}\\
Fig. 2. \textit{Pair--of--pants object that defines product on $HF_*(o_M,\nu^*N:H)$}
\end{center}
\vskip5mm

When $N=M$ we obtain the product defined in~\cite{O2}. It was proven in~\cite{KMS} that the induced operation on $H_*(M)$ is exactly the intersection product. We can describe the operation on singular homology as the operation given by composing the appropriate action and the inclusion morphism (see Section~\ref{products}).

The product $\ast$ can be used in order to prove a triangle inequality for conormal spectral invariant. Our inequality is a generalization of the one made by Monzner, Vichery and Zapolsky in \cite{MNZ}.
\begin{prop}\label{inv_ineq}
Let us take two compactly supported Hamiltonians $H,H'$ and $\alpha,\beta\in H_*(N)$ such that $\alpha\cdot\beta\neq0$. Then
$$
l(\alpha\cdot\beta;o_M,\nu^*N:H\sharp H')\le l(\alpha;o_M,\nu^*N:H)+l(\beta;o_M,\nu^*N:H').
$$
\end{prop}
\bigskip
\indent If we now take a pair--of--pants with different type of boundary conditions (see figure 3) we can prove that conormal spectral invariants are bounded for every non--zero singular homology class. The idea of this property came from a Humili\`{e}re, Leclercq and Seyfaddini's paper (see ~\cite{HLS}).
\begin{prop}\label{bound}
For every $\alpha\in H_*(N)\setminus\{0\}$ it holds
$$l(\alpha;o_M,\nu^*M:H)\le l([M];o_M,o_M:H),
$$
where $[M]\in H_*(M)$ is the fundamental class.
\end{prop}
\vskip2mm
\begin{center}\includegraphics[width=6 cm,height=3cm]{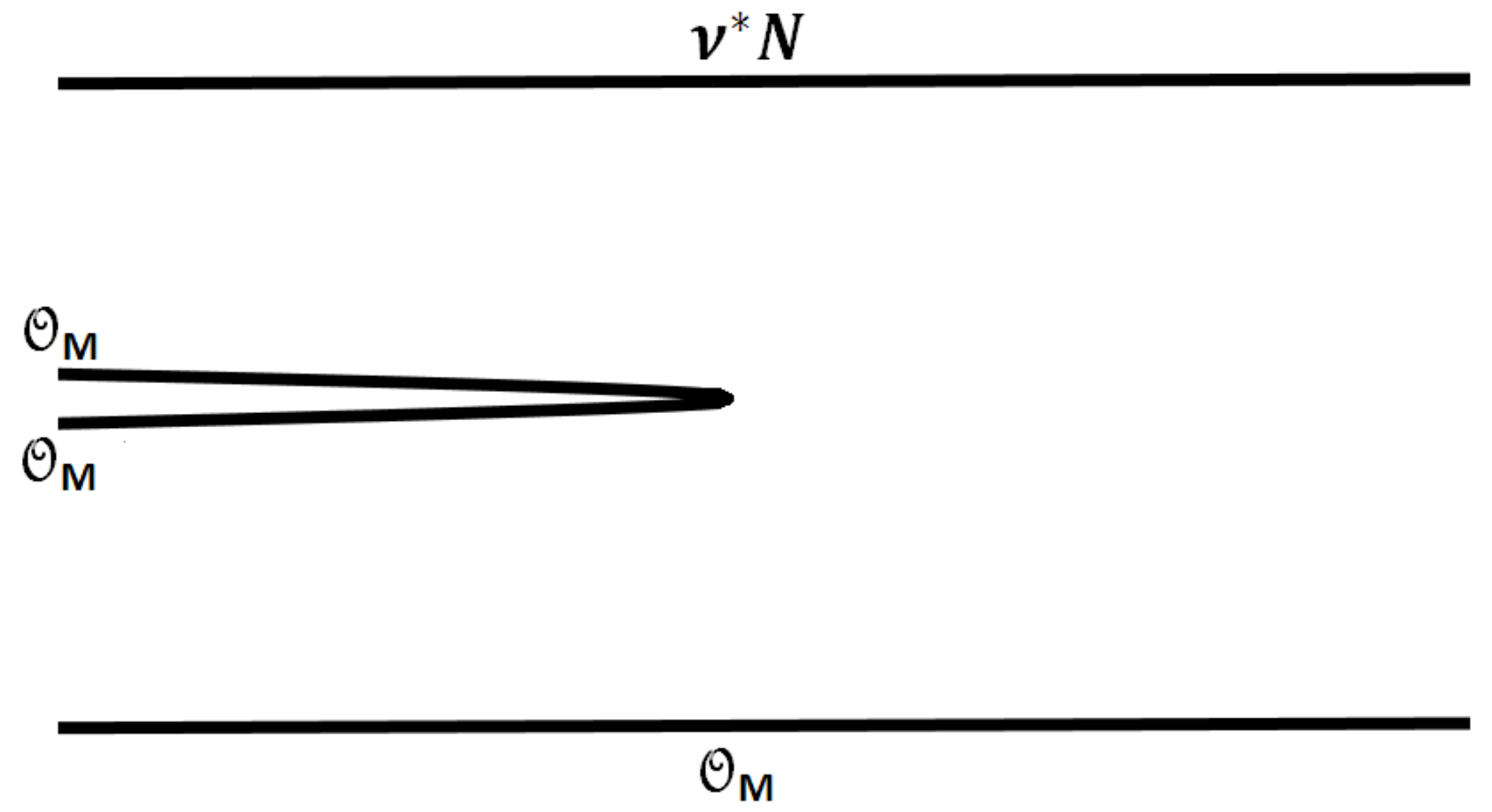}\\
Fig. 3. \textit{Pair--of--pants object that gives boundness of conormal spectral invariants}
\end{center}
\vskip5mm

This paper is organized as follows. In Section~\ref{diskovi i moduli} we define diverse moduli spaces and prove some of their properties. In Section~\ref{izomorfizam} we present the construction of PSS--type homomorphisms and we prove Theorem~\ref{main_thm}. Section~\ref{dijagram} contains the proof of Theorem~\ref{diagram_commutes}. In the last section we give the construction of a product in homology and prove Proposition~\ref{inv_ineq} and Proposition~\ref{bound}.\\

\smallskip
\textbf{Acknowledgments.} The author thanks Jelena Kati\'{c}, Darko Milinkovi\'{c} and Katrin Wehrheim for useful discussions during the preparation of this paper. The author also thanks the anonymous referee for many valuable suggestions and corrections.
\smallskip
\section{\textbf{Holomorphic discs, gradient trajectories and moduli spaces}}\label{diskovi i moduli}

We start with a construction of mixed--type object space that we use for the definition of $\Psi$ and $\Phi$. Let $p$ be a critical point of a Morse function $f$. Morse homology $HM_k(f)$ is graded by Morse index $k=m_f(p)$ of critical points.\\
\indent To each element of $CF_*(H)$ we can assign a solution of Hamiltonian equation
\begin{equation}\label{eq:HamSys}
\left\{
\begin{array}{ll}
\dot{x}=X_H(x),\\
x(0)\in o_M,\,\,x(1)\in \nu^*N.\\
\end{array}
\right.
\end{equation}
For a solution $x$ of~(\ref{eq:HamSys}) there exists a canonically assigned Maslov index $$\mu_N:CF_*(H)\rightarrow\frac{1}{2}\mathbb{Z},$$ see~\cite{O1,RS1,RS2} for details. Floer homology $HF_k(H)$ is graded by $k=\mu_N(x)+\frac{1}{2}\dim N$.

Let ${\mathcal M} (p,f;x,H)$ be the space of pairs of maps
$$\gamma:(-\infty, 0]\to N, \quad u:[0,+\infty)\times [0,1]\to T^*M,$$
that satisfy
\begin{equation}\label{eq:grad1}
\left\{
\begin{array}{ll}
\frac{d\gamma}{ds}=-\nabla f(\gamma(s)), \\
\frac{\partial u}{\partial s}+J(\frac{\partial u}{\partial t}-X_{\rho_R^+ H}(u))= 0, \\
u(s,0)\in o_M, u(s,1)\in\nu^*N, u(0,t)\in o_M, s\geq0, t\in[0,1], \\
\gamma(-\infty)=p,
u(+\infty,t)=x(t), \\
\gamma(0)=u(0,1),\\
\end{array}
\right.
\end{equation}
where $R$ is a positive fixed number and $\rho_R^+:[0,+\infty)\rightarrow{\mathbb R}$ is a smooth function such that
$$
\rho_R^+(s)=\begin{cases} 1, & s\ge R+1, \\ 0, & s\le R. \end{cases}
$$

Let ${\mathcal M} (x,H;p,f)$ be the space of pairs of maps
$$\gamma:[0,+\infty)\to N, \quad u:(-\infty,0]\times [0,1]\to T^*M,$$
that satisfy
\begin{equation}\label{eq:grad1}
\left\{
\begin{array}{ll}
\frac{d\gamma}{ds}=-\nabla f(\gamma(s)), \\
\frac{\partial u}{\partial s}+J(\frac{\partial u}{\partial t}-X_{\rho_R^- H}(u))= 0, \\
u(s,0)\in o_M, u(s,1)\in\nu^*N, u(0,t)\in o_M, s\leq0, t\in[0,1], \\
\gamma(+\infty)=p,
u(-\infty,t)=x(t), \\
\gamma(0)=u(0,1),\\
\end{array}
\right.
\end{equation}
where $\rho_R^-:(-\infty,0]\rightarrow{\mathbb R}$ is a smooth function such that
$$
\rho_R^-(s)=\begin{cases} 1, & s\le -R-1, \\ 0, & s\ge -R. \end{cases}
$$
\vskip5mm
\begin{center}\includegraphics[width=6.2cm,height=4cm]{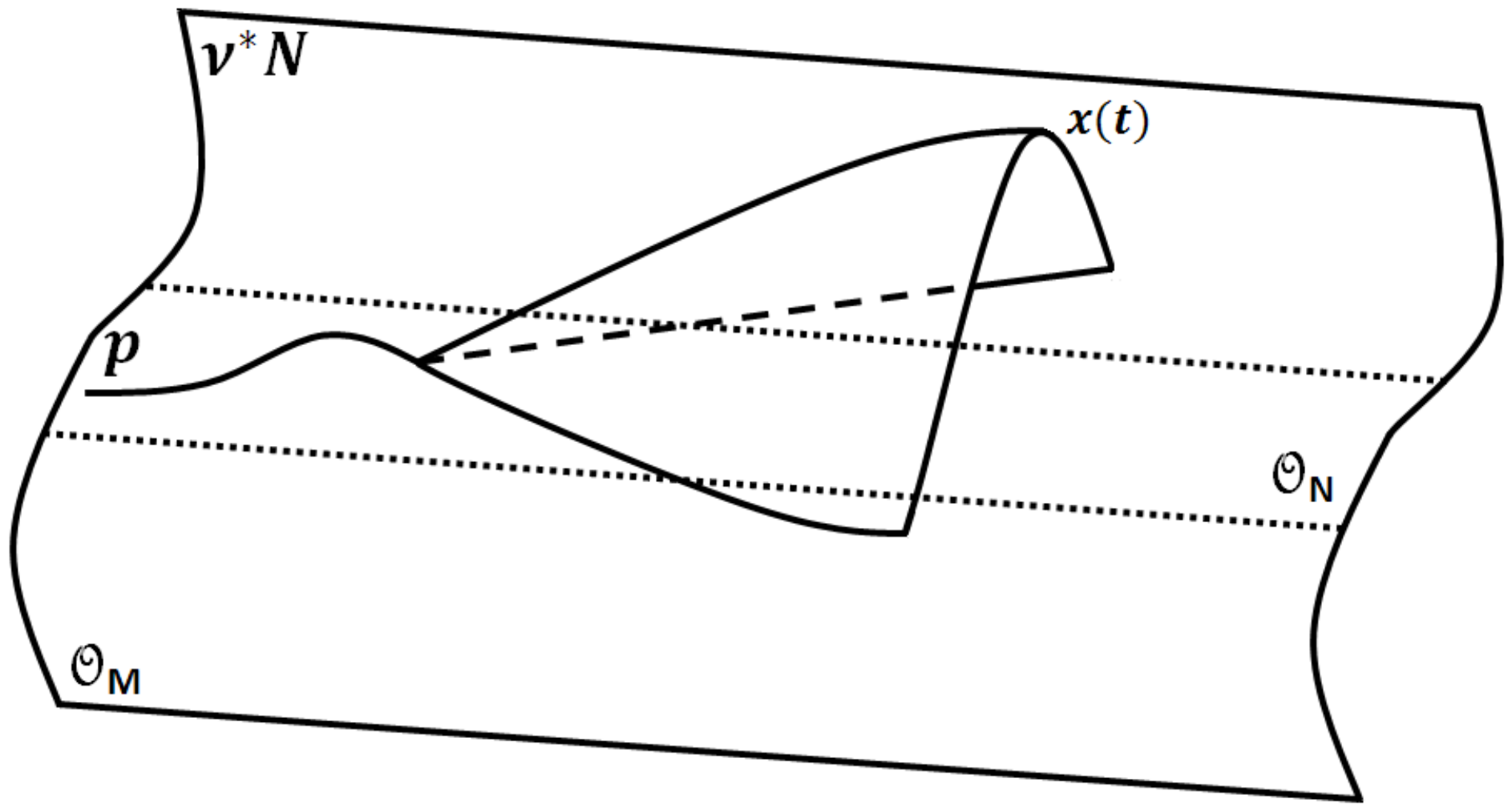}\hskip2mm\includegraphics[width=6.2cm,height=4cm]{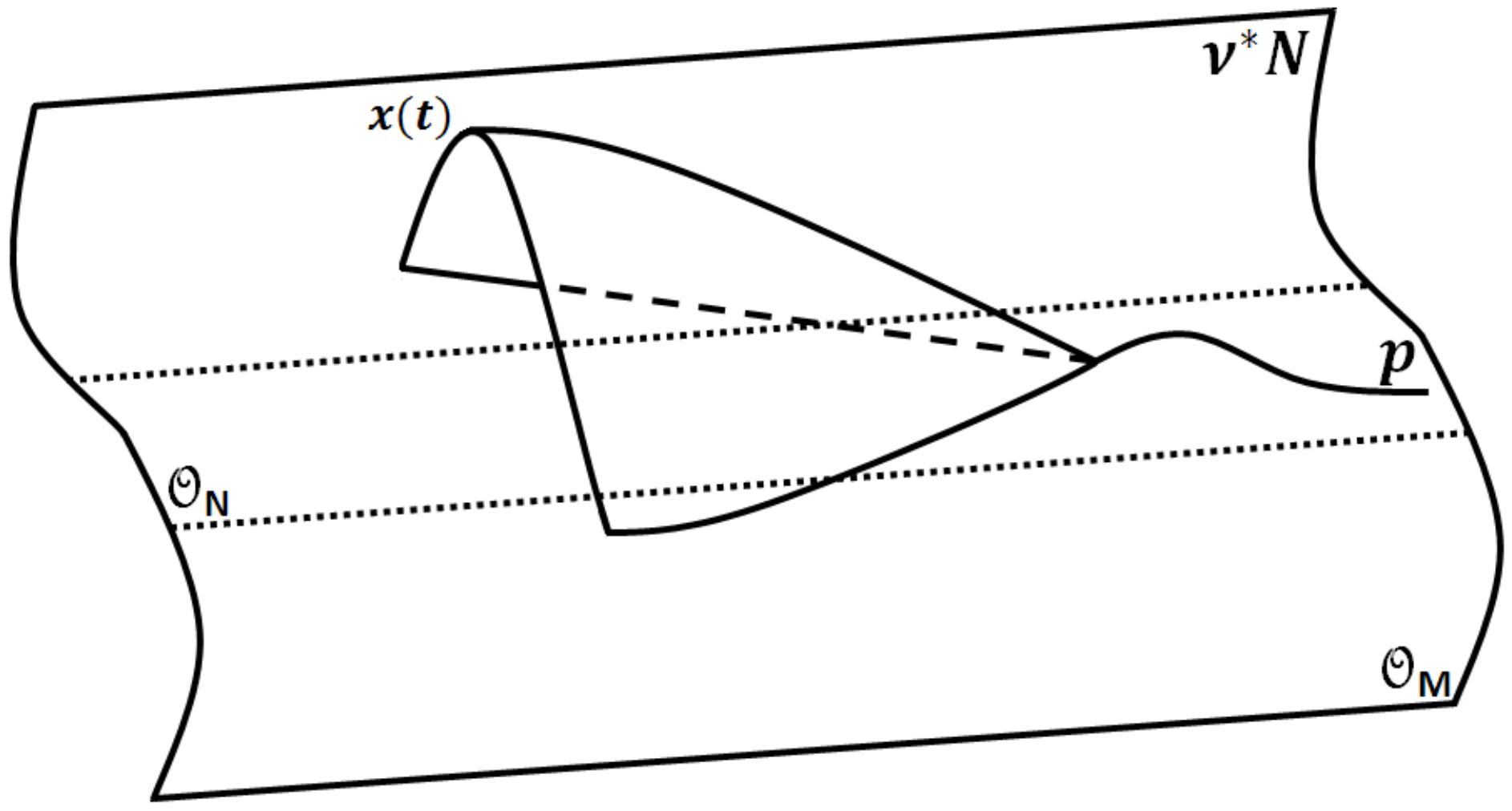}\\
Fig. 4. ${\mathcal M}(p,f;x,H)$ \textit{and} ${\mathcal M}(x,H;p,f)$
\end{center}
\vskip5mm
\begin{prop}\label{dim}
For a generic Morse function $f$ and a generic compactly supported Hamiltonian $H$ the set ${\mathcal M}(p,f;x,H)$ is a smooth manifold
of dimension $m_f(p)-(\mu_N(x)+\frac{1}{2}\dim N)$ and ${\mathcal M}(x,H;p,f)$ is
a smooth manifold of dimension $\mu_N(x)+\frac{1}{2}\dim N-m_f(p)$.
\end{prop}
\smallskip

\noindent {\it Proof:} Let $W^u(p,f)$ be the unstable manifold associated to a critical point
$p$ of a Morse function $f$. We know that $\dim W^u(p,f)=m_f(p)$ (see~\cite{M}). Let $W^s(x,H)$ be the set of solutions of

\begin{equation}\label{halfStrip}
\left\{
\begin{array}{ll}
u:[0,+\infty)\times[0,1]\rightarrow T^*M,\\
\frac{\partial u}{\partial s}+J(\frac{\partial u}{\partial t}-X_{\rho_R^+H}(u))= 0, \\
u(s,0)\in o_M, u(s,1)\in\nu^*N, u(0,t)\in o_M, s\geq 0, t\in[0,1], \\
u(+\infty,t)=x(t).
\end{array}
\right.
\end{equation}

Now we compute the dimension of $W^s(x,H)$. It goes similarly as in~\cite{AS} (Theorem 5.24 and Theorem 5.25), but for the sake of completeness, we give the main steps adopted to our situation. \\
Let us denote $V={\mathbb R}^m\times\{0\}\subset{\mathbb C}^m$ and let $W=W_0\times{\mathbb R}^m\subset{\mathbb C}^m$ be a linear subspace (in our case $W_0$ will be a local model for $N$). For convenience, we use the notation
$$\Sigma^+=\big\{z\in{\mathbb C}\,|\,{\mathrm Re}{z}\ge0,\,0\le{\mathrm Im}{z}\le1\big\},
$$
instead of $[0,+\infty)\times[0,1]$ and $s+it$ instead of $(s,t)\in[0,+\infty)\times[0,1]$. Let $X^{1,p}_W(\Sigma^+,{\mathbb C}^m)$ be the completion of the space of maps with bounded support $$
\begin{aligned}
&u\in C^\infty_c(\Sigma^+,{\mathbb C}^m),\\
&u(it)\in V, t\in[0,1],\hskip5mm\\
&(u(s),\overline{u}(s+i))\in\nu^*W, s\ge0,
\end{aligned}
$$
with the respect to the norm
$$
\big\|u\big\|_{X^{1,p}(\Sigma^+)}=\big\|u\big\|^p_{L^p(\Sigma^+)}+\big\|Du\big\|^p_{L^p(\Sigma^+)}.
$$
The space $X^p(\Sigma^+,{\mathbb C}^m)$ is the space of locally integrable ${\mathbb C}^m$--valued maps on $\Sigma^+$ whose $\big\|\cdot\big\|_{L^p(\Sigma^+)}$ norm is finite. Let $A\in C^0([0,+\infty)\times[0,1],L({\mathbb R}^{2n},{\mathbb R}^{2n}))$ be such that $A(+\infty,t)\in\Sym(2n,{\mathbb R})$ for every $t\in[0,1]$ (i.e. $A(+\infty,t)$ is symmetric). Denote by $\Phi^+:[0,1]\to\Sp(2n)$ the solutions of the linear Hamiltonian system
$$
\frac{d}{dt}\Phi^+(t)=iA(+\infty,t)\Phi^+(t),\hskip3mm \Phi^+(0)={\mathbb I}.
$$
Then, (see~\cite{AS}) for every $p\in(1,+\infty)$ the operator
$$
\begin{aligned}
&\overline{\partial}_A:X^{1,p}_W(\Sigma^+,{\mathbb C}^m)\to X^p(\Sigma^+,{\mathbb C}^m),\\
&\overline{\partial}_A u=\partial_s u+i\partial_t u+Au,
\end{aligned}
$$
is bounded and Fredholm of index
$$
\ind\overline{\partial}_A=\frac{1}{2}\dim W_0-\mu(\verb"graph"{\,\Phi^+C},\nu^*W).
$$
Here, $C$ denotes the anti--symplectic involution of $T^*{\mathbb R}^m$ which maps $(q,p)$ into $(q,-p)$ and $\mu$ denotes the relative Maslov index of two paths of Lagrangian subspaces of $T^*{\mathbb R}^{2m}$ (defined in~\cite{RS1, RS2}). We can see $W^s(x,H)$ as the set of zeroes of a smooth section of a suitable Banach bundle. The fiberwise derivative of such section at $u\in W^s(x,H)$ is conjugated to a linear operator $\overline{\partial}_A$. It follows that dimension of $W^s(x,H)$ is equal to the Fredholm index of an operator $\overline{\partial}_A$ with $V={\mathbb R}^{m}$, $m=\dim M$ and $W$ a local model for $N\times M$. Thus,
$$\dim W^s(x,H)=\frac{1}{2}\dim N-\mu_N(x).$$
We used the definition of Maslov index
$$
\mu_N(x)=\mu(B_{\Phi}({\mathbb R}^{m}),V^{\Phi}),
$$
where $\Phi:x^*T(T^*M)\to[0,1]\times{\mathbb C}^m$ is any trivialization and
$$
\begin{aligned}
&V^{\Phi}=\Phi(T_{x(1)}\nu^*N),\\
&B_\Phi(t)=\Phi\circ T\phi_H^t\circ\Phi^{-1}.
\end{aligned}
$$
For a generic choice of parameters the evaluation map $$Ev:W^u(p,f)\times W^s(x,H)\rightarrow N\times N,\,\,\,Ev(\gamma,u)=(\gamma(0),u(0,1)),$$ is transversal to the diagonal, thus ${\mathcal M}(p,f;x,H)=Ev^{-1}(\triangle)$ is a smooth manifold of dimension $m_f(p)+\frac{1}{2}\dim N-\mu_N(x)-(2\dim N-\dim N)=m_f(p)-\frac{1}{2}\dim N-\mu_N(x)$. The proof for ${\mathcal M}(x,H;p,f)$ is similar.
\qed

\begin{rem}\label{napomena}Referee suggested that previous construction (of half--strips in~(\ref{halfStrip})) can be viewed in more elegant way as follows. Observing Floer equation on $\dot{\mathbb{D}^2}=\mathbb{D}^2\setminus\{1\}\subset{\mathbb{C}}$ where the bottom half of $\partial\dot{\mathbb{D}}^2$ goes to $o_M$, and the top half of $\partial\dot{\mathbb{D}}^2$ goes to $\nu^*N$. The Floer equation can be written relative to the standard conformal coordinates on $\mathbb{D}^2\setminus\{\pm1\}$ given by the biholomorphism
$$
\mathbb{R}\times[0,1]\to\mathbb{D}^2\setminus\{\pm1\},\,\,\,z=s+it\mapsto\frac{e^{\pi z}-i}{e^{\pi z}+i}.
$$
Relative to the coordinates $(s,t)$ on $\mathbb{D}^2\setminus\{\pm1\}$ induced by it the Floer equation is the usual one where the Hamiltonian term is cut off for $s$ large enough. There is no issue at $-1$ since there the equation is just the Cauchy–-Riemann equation, which is independent of the conformal coordinates chosen.

Oh used half--strips with boundary on the zero section in order to find the dimension of moduli space of pair--of--pants with boundary on the zero section (see Appendix in~\cite{O2}, Theorem A.1). We used similar half--strips with switching boundary conditions as Abbondandolo and Schwarz (see~\cite{AS}, Corollary 5.30).
\end{rem}



We need some additional properties on manifolds ${\mathcal M}(p,f;x,H)$ and ${\mathcal M}(x,H;p,f)$. The set of solutions of~(\ref{eq:1}) is denoted by ${\mathcal M}(x,y;H)$ and ${\mathcal M}(p,q;f)$ denotes the set of solutions of~(\ref{eq:3}) (modulo ${\mathbb R}$--action).

\begin{prop}\label{boundary1}
Let $f$ be a generic Morse function and $H$ a generic compactly supported Hamiltonian. If $m_f(p)=\mu_N(x)+\frac{1}{2}\dim N$ then ${\mathcal M}(p,f;x,H)$ is a finite set. If $m_f(p)=\mu_N(x)+\frac{1}{2}\dim N+1$ then ${\mathcal M}(p,f;x,H)$ is one--dimensional manifold with topological boundary

$$\begin{aligned}
\partial{\mathcal M}(p,f;x,H)=
&\bigcup_{m_f(q)=m_f(p)-1}\mathcal{M}(p,q;f)\times{\mathcal M}(q,f;x,H)\\
\cup\,\,&\bigcup_{\mu_N(y)=\mu_N(x)+1}\mathcal{M}(p,f;y,H)\times{\mathcal M}(y,x;H).
\end{aligned}$$

\end{prop}
\smallskip
\noindent {\it Proof:}
Let $(\gamma_n,u_n)$ be a sequence in ${\mathcal M}(p,f;x,H)$ that has no $W^{1,2}-$convergent subsequence. Since $N$ is compact $\gamma_n(t)$ is bounded for every $t$. The sequence $\gamma_n$ is equicontinuous because
$$
\begin{aligned}
d(\gamma_n(t_1),\gamma_n(t_2))&\le\int_{t_1}^{t_2}\|\dot{\gamma}(s)\|\,ds\\
&\le\sqrt{t_2-t_1}\sqrt{\int_{t_1}^{t_2}\|\dot{\gamma}(s)\|^2\,ds}\\
&=\sqrt{t_2-t_1}\sqrt{\int_{t_1}^{t_2}\frac{\partial}{\partial s}f(\gamma_n(s))\,ds}\\
&\le\sqrt{t_2-t_1}\sqrt{\max_{x\in N}f(x)-f(\gamma_n(-\infty))}\\
&=\sqrt{t_2-t_1}\sqrt{\max_{x\in N}f(x)-f(p)}.
\end{aligned}
$$
It follows from the Arzel\`{a}--Ascoli theorem that $\gamma_n$ has a subsequence which converges uniformly on compact sets. Since this sequence $\gamma_n$ is a solution of an equation
$$
\dot{\gamma}_n=-\nabla f(\gamma_n),
$$
and function $f$ is smooth, $\gamma_n$ converges with all its derivatives on compact subsets of $(-\infty,0]$.

The energy of $u_n$,
$$
E(u_n)=\int_0^{+\infty}\int_0^1\bigg\|\frac{\partial u_n}{\partial s}\bigg\|^2+\bigg\|\frac{\partial u_n}{\partial t}-X_{\rho_R^+H}(u_n)\bigg\|^2\,dt\,ds,
$$
is uniformly bounded.
$$
E(u_n)=\bigg(\int_0^{R+1}\int_0^1+\int_{R+1}^{+\infty}\int_0^1\bigg)\bigg\|\frac{\partial u_n}{\partial s}\bigg\|^2+\bigg\|\frac{\partial u_n}{\partial t}-X_{\rho_R^+H}(u_n)\bigg\|^2\,dt\,ds.
$$
Uniform bound of the first integral follows from an estimate
$$
\|u\|_{W^{1,2}([0,1]\times[0,R+1])}\le c_0\|u\|_{L^2(U)}+c_1\|\bar{\partial}u\|_{L^2(U)},
$$
and from $C^0$--boundness of a sequence $u_n$ (see \cite{AS} for details). Here, $U$ is some open subset of finite measure that contains $[0,1]\times[0,R+1]$.
The second integral is uniformly bounded because
$$
\begin{aligned}
\int_{R+1}^{+\infty}\int_0^1\bigg\|\frac{\partial u_n}{\partial s}\bigg\|^2+&\bigg\|\frac{\partial u_n}{\partial t}-X_{\rho_R^+H}(u_n)\bigg\|^2\,dt\,ds\\
&=\int_{R+1}^{+\infty}\int_0^1\bigg\|\frac{\partial u_n}{\partial s}\bigg\|^2+\bigg\|\frac{\partial u_n}{\partial t}-X_H(u_n)\bigg\|^2\,dt\,ds\\
&\le\int_{0}^{+\infty}\int_0^1\bigg\|\frac{\partial u_n}{\partial s}\bigg\|^2+\bigg\|\frac{\partial u_n}{\partial t}-X_H(u_n)\bigg\|^2\,dt\,ds\\
&={\mathcal A}_H(x(\cdot))-{\mathcal A}_H(u_n(0,\cdot))\\
&\le{\mathcal A}_H(x(\cdot))-\min{H}.
\end{aligned}
$$
We have a sequence $u_n$ whose energy is uniformly bounded. From Gromov compactness (see~\cite{G}) it follows that $u_n$ has a subsequence that converges together with all derivatives on compact subsets of $([0,+\infty)\times[0,1])\setminus\{z_1,...z_m\}$. Bubbles can occur at $z_i$ if it is an interior point of $[0,+\infty)\times[0,1]$. It is also possible that a bubble appears at the boundary point $z_k$ as holomorphic disc with the boundary conditions on zero section and conormal bundle. But in our case neither holomorphic spheres nor discs appear. If $v:S^2\to T^*M$ is a holomorphic sphere then
$$
\int_{S^2}\|dv\|^2=\int_{S^2}v^*\omega=\int_{\partial S^2}v^*\lambda=0.
$$
If $v:[0,+\infty)\times[0,1]\to T^*M$ is a holomorphic disc then
$$
\int_{[0,+\infty)\times[0,1]}\|dv\|^2=\int_{[0,+\infty)\times[0,1]}v^*\omega=\int_{\partial([0,+\infty)\times[0,1])}v^*\lambda=0,
$$ since $\lambda=0$ on $o_M$ and $\nu^*N$.\\
\indent So, $(\gamma_n,u_n)$ has a subsequence which converges with all its derivatives uniformly on compact sets. From $C^{\infty}_{loc}$ convergence it follows $W^{1,2}$ convergence. Thus, $(\gamma_n,u_n)$ has a subsequence that converges to some element of ${\mathcal M}(p^m,f;x^0,H)$. Similarly as in~\cite{F2,K,MS,Sa,Sc1} we conclude that the only loss of compactness is a "trajectory breaking" in the following way
\begin{equation}\label{big_union}
\begin{aligned}
\bigcup\,&{\mathcal M}(p,p^1;f)\times...\times{\mathcal M}(p^{m-1},p^m;f)\times{\mathcal M}(p^m,f;x^0,H)\\
\times&{\mathcal M}(x^0,x^1;H)\times...\times{\mathcal M}(x^{l-1},x;H).
\end{aligned}
\end{equation}
Here, $p,p^1,...,p^m$ are critical points of $f$ and $x^0,...,x^{l-1},x$ are Hamiltonian paths with decreasing Morse and Maslov indices such that $m_f(p^m)\ge\mu_N(x^0)+\frac{1}{2}\dim N$. Therefore, we have that a boundary $\partial{\mathcal M}(p,f;x,H)$ is a subset of an union in~(\ref{big_union}). The other inclusion follows from standard gluing arguments.\\
\indent If $m_f(p)=\mu_N(x)+\frac{1}{2}\dim N$ then ${\mathcal M}(p,f;x,H)$ is a compact, zero--dimensional manifold, so ${\mathcal M}(p,f;x,H)$ has a finite number of elements.\\
\indent If $m_f(p)=\mu_N(x)+\frac{1}{2}\dim N+1$ then the boundary of ${\mathcal M}(p,f;x,H)$ can contain an element of a set ${\mathcal M}(p,q;f)\times{\mathcal M}(q,f;x,H)$ for some $q\in\Crit(f)$ such that $m_f(q)=m_f(p)-1$ or an element of a set ${\mathcal M}(p,f;y,H)\times{\mathcal M}(y,x;H)$ for some Hamiltonian orbit $y$, such that $\mu_N(y)=\mu_N(x)+1$.
\qed
\\

We have a similar proposition for ${\mathcal M}(x,H;p,f)$.
\begin{prop}\label{boundary2}
Let $f$ be a generic Morse function and $H$ a generic compactly supported Hamiltonian. If $m_f(p)=\mu_N(x)+\frac{1}{2}\dim N$ then ${\mathcal M}(x,H;p,f)$ is a finite set. If $m_f(p)=\mu_N(x)+\frac{1}{2}\dim N-1$ then ${\mathcal M}(x,H;p,f)$ is one--dimensional manifold with topological boundary
$$\begin{aligned}
\partial{\mathcal M}(x,H;p,f)=
&\bigcup_{m_f(q)=m_f(p)+1}\mathcal{M}(x,H;q,f)\times{\mathcal M}(q,p;f)\\
\cup\,\,&\bigcup_{\mu_N(y)=\mu_N(x)-1}\mathcal{M}(x,y;H)\times{\mathcal M}(y,H;p,f).
\end{aligned}$$\qed
\end{prop}
Now, we define some auxiliary manifolds that we use to prove that the composition $\Phi\circ\Psi$ is the identity (see Theorem~\ref{main_thm}). Let $R>0$ be a fixed number. For $p,q\in\Crit(f)$ define
$$
{\mathcal M}_R(p,q,f;H)=
\left\{ (\gamma_{-},\gamma_{+},u) \left|
\begin{array}{ll}
\gamma_{-}:(-\infty,0]\rightarrow N, \;
\gamma_{+}:[0,+\infty)\rightarrow N, \\
u:{\mathbb R}\times [0,1]\rightarrow T^*M, \\
\frac{d\gamma_{\pm}}{ds}=-\nabla f(\gamma_{\pm}), \\
\frac{\partial u}{\partial s}+J(\frac{\partial u}{\partial t}-X_{\sigma_RH}(u))= 0, \\
\gamma_-(-\infty)=p, \; \gamma_+(+\infty)=q, \\
u(s,0)\in o_M, u(s,1)\in\nu^*N,s\in{\mathbb R},\\
u(-\infty,t),u(+\infty,t)\in o_M, t\in[0,1],\\
u(\pm\infty,t)=\gamma_{\pm}(0)
\end{array}
\right.
\right\},
$$
where $\sigma_R:{\mathbb R}\rightarrow[0,1]$ is a smooth function such that
$$
\sigma_R(s)=\begin{cases} 1, & |s|\le R, \\ 0, & |s|\ge R+1, \end{cases}
$$ and $$\overline{{\mathcal M}}(p,q,f;H)=\big\{(R,\gamma_-,\gamma_+,u)\,|\,(\gamma_-,\gamma_+,u)\in{\mathcal M}_R(p,q,f;H),\,R>R_0\big\},$$ (see figure 5). For a generic choice of parameters, the set $\overline{{\mathcal M}}(p,q,f;H)$ is an one--dimensional manifold if $m_f(p)=m_f(q)$ and a zero--dimensional manifold if $m_f(p)=m_f(q)-1$.
\vskip4mm
\begin{center}\includegraphics[width=9cm,height=4.5cm]{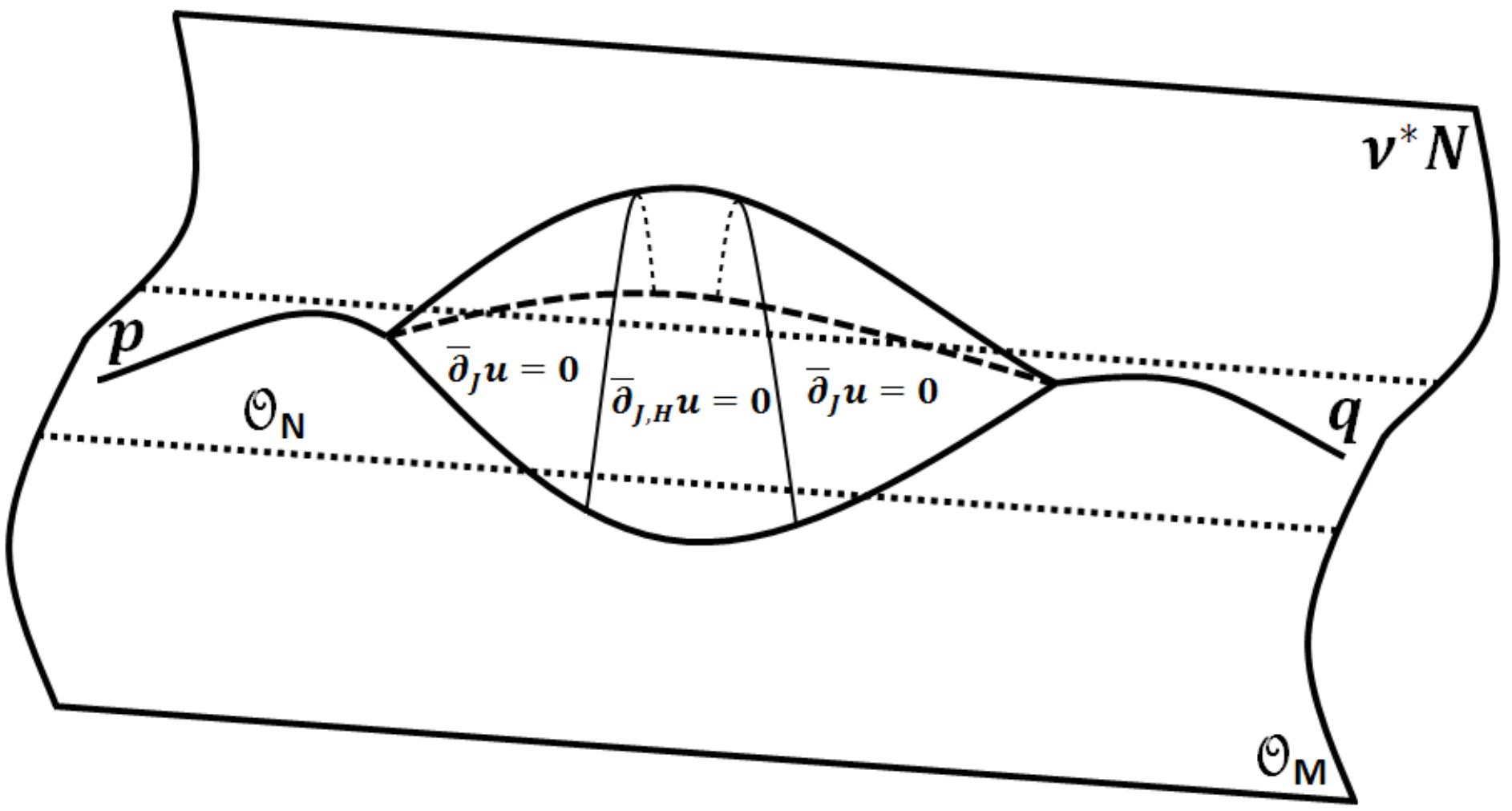}\\
Fig. 5. ${\mathcal M}_R(p,q,f;H)$
\end{center}
\vskip5mm
Knowing the definitions of a broken gradient trajectory and a weak convergence of gradient trajectories (see~\cite{Sc1}) we can define a broken holomorphic strip and a weak convergence of holomorphic strips (see~\cite{Sa}).
\begin{defn}\label{broken_strip} A broken (perturbed) holomorphic strip $v$ is a pair $(v_1,v_2)$ of (perturbed) holomorphic strips such that $v_1(+\infty,t)=v_2(-\infty,t)$. A sequence of perturbed holomorphic strips $u_n:{\mathbb R}\times[0,1]\to T^*M$ is said to converge weakly to a broken trajectory $v$ if there exists a sequence of translations $\varphi_n^i:{\mathbb R}\times[0,1]\to{\mathbb R}\times[0,1]$, $i=1,2$, such that $u_n\circ\varphi_n^i$ converges to $v_i$ uniformly with all derivatives on compact subset of ${\mathbb R}\times[0,1]$. We say that an element of mixed type $(\gamma,u)$ is a broken element if $\gamma$ is a broken trajectory or $u$ is a broken holomorphic strip.
\end{defn}
Next proposition gives us a boundary of an one--dimensional manifold $\overline{{\mathcal M}}(p,q,f;H)$.
\begin{prop}\label{boundary3}
Let $p,q\in CM_k(f)$. Then the topological boundary of $\overline{{\mathcal M}}(p,q,f;H)$ can be identified with
$$
\begin{aligned}
\partial\overline{\mathcal M}(p,q,f;H)={\mathcal M}_{R_0}(p,q,f;H)\,\,\,&\cup\,\,\,\,\,\,\,\,\bigcup_{m_f(r)=k-1}{\mathcal M}(p,r;f)\times\overline{\mathcal M}(r,q,f;H)\\
&\cup\,\,\,\,\,\,\,\,\bigcup_{m_f(r)=k+1}\overline{\mathcal M}(p,r,f;H)\times{\mathcal M}(r,q;f)\\
&\cup\bigcup_{\mu_N(x)+\dim N/2=k}{\mathcal M}(p,f;x,H)\times{\mathcal M}(x,H;q,f).
\end{aligned}
$$
\end{prop}
\smallskip
\noindent {\it Proof:} Let us take a sequence $(R_n,\gamma_-^n,\gamma_+^n,u_n)$ in $\overline{{\mathcal M}}(p,q,f;H)$. Then, this sequence either $W^{1,2}-$converges to an element of the same moduli space or one of the following four statements holds:
\begin{itemize}
\item[(1)] There is a subsequence such that $R_{n_k}\rightarrow R_0$ and $(\gamma_-^{n_k},\gamma_+^{n_k},u_{n_k})$ converges to $(\gamma_-,\gamma_+,u)\in{\mathcal M}_{R_0}(p,q,f;H)$.
\item[(2)] There is a subsequence of $(R_n,\gamma_-^n,\gamma_+^n,u_n)$ that converges to a broken trajectory in ${\mathcal M}(p,r;f)\times\overline{\mathcal M}(r,q,f;H)$. Subsequence $(\gamma_+^{n_k},u_{n_k})$ converges in $W^{1,2}$ topology and $\gamma_-^{n_k}$ converges weakly.
\item[(3)] There is a subsequence that converges to a broken trajectory in $\overline{\mathcal M}(p,r,f;H)\times{\mathcal M}(r,q;f)$, similarly to (2).
\item[(4)] There is a subsequence such that $R_{n_k}\rightarrow+\infty$ and $(\gamma_-^{n_k},\gamma_+^{n_k},u_{n_k})$ converges weakly to a broken element of ${\mathcal M}(p,f;x,H)\times{\mathcal M}(x,H;q,f)$.
\end{itemize}
If $R_n$ is bounded then we can find compact $K$ such that $\{R_n\}\subset K$. The family $\rho_R$ can be chosen to depend continuously on $R$, so all estimates in Proposition \ref{boundary1} hold uniformly on $R\in K$. In a similar way to Proposition \ref{boundary1} we conclude that $(\gamma_-^n,\gamma_+^n,u_n)$ has a subsequence that converges locally uniformly. So, if $(R_n,\gamma_-^n,\gamma_+^n,u_n)$ does not converge to an element of $\overline{{\mathcal M}}(p,q,f;H)$, then $R_n\to R_0$ or $R_n\to R>R_0$ ($R_n$ denotes the subsequence, as well). If the first case $(\gamma_-^n,\gamma_+^n,u_n)$ converges in $W^{1,2}$ topology and in the second one $(\gamma_-^n,\gamma_+^n,u_n)$ converges to a broken trajectory. Since dimension of $\overline{{\mathcal M}}(p,q,f;H)$ is one it can break only once. The breaking can happen on trajectories $\gamma_-^n$ or $\gamma_+^n$ and not on the disc. Sequence $u_n$ cannot converge to a broken disc because the non--holomorphic part of the domain is compact and there $u_n$ converges. If it breaks on the holomorphic part we obtain a solution of a system
$$
\left\{
\begin{array}{ll}
v:{\mathbb R}\times[0,1]\to T^*M,\\
\frac{\partial v}{\partial s}+J\frac{\partial v}{\partial t}=0,\\
v({\mathbb R}\times\{0\})\subset O_M,\,\,v({\mathbb R}\times\{1\})\subset \nu^*N.
\end{array}
\right.
$$
We already saw that all such solutions are constant, so $u_n$ cannot break on the holomorphic part neither. In this way we covered the first three cases. The fourth case arises if $R_n$ is not bounded sequence. We can find a subsequence $R_n\to+\infty$. Then discs
$$
u_n^-(s,t):=u_n(s-R_n-R_0-1,t),\,\,u_n^+(s,t):=u_n(s+R_n+R_0+1,t),$$
converge locally uniformly with all derivatives to some $u^-$ and $u^+$. These discs are solutions of the system
$$
\left\{
\begin{array}{ll}
\frac{\partial u^{\pm}}{\partial s}+J(\frac{\partial u^{\pm}}{\partial t}-X_{\rho^\pm_{R_0}}(u^\pm))=0,\\
u^\pm({\mathbb R}\times\{0\})\subset O_M,\,\,u^\pm({\mathbb R}\times\{1\})\subset\nu^*N,\\
u^\pm(\mp\infty,t)=x(t),\\
u^\pm(\pm\infty,t)=\gamma_\pm(0).
\end{array}
\right.
$$
Sequences $\gamma_\pm^n$ cannot break because of dimensional reason so they converge to some trajectories $\gamma_\pm$.
\\
Conversely, for each broken trajectory of some of these types:
\begin{itemize}
\item $(\gamma,\gamma_-,\gamma_+,u)\in{\mathcal M}(p,r;f)\times\overline{\mathcal M}(r,q,f;H),$
\item $(\gamma_-,\gamma_+,u,\gamma)\in\overline{\mathcal M}(p,r,f;H)\times{\mathcal M}(r,q;f),$
\item $(\gamma_1,u_1,\gamma_2,u_2)\in{\mathcal M}(p,f;x,H)\times{\mathcal M}(x,H;q,f)$,
\end{itemize}
there is a sequence in $\overline{{\mathcal M}}(p,q,f;H)$ that converges weakly to a corresponding broken trajectory. The proof is based on the implicit function theorem and pre--gluing and gluing techniques.\qed\\

We continue with the construction of the auxiliary manifold, again with the variable domain, that now connects Hamiltonian orbits. Let $\varepsilon>0$ be a fixed number. Consider
$$
{\mathcal M}_{\varepsilon}(x,y,H;f)=\left\{ (u_{-},u_{+},\gamma) \left|
\begin{array}{ll}
u_{-}:(-\infty,0]\times[0,1]\rightarrow T^*M,\\
u_{+}:[0,+\infty)\times[0,1]\rightarrow T^*M, \\
\gamma:[-\varepsilon,\varepsilon]\rightarrow N, \\
\frac{\partial u_{\pm}}{\partial s}+J(\frac{\partial u_{\pm}}{\partial t}-X_{\rho_R^{\pm}H}(u_{\pm}))= 0, \\
\frac{d\gamma}{ds}=-\nabla f(\gamma), \\
u_{\pm}(\pm s,0)\in o_M, u_{\pm}(\pm s,1)\in\nu^*N,s\ge0,\\
u_{\pm}(0,t)\in o_M, t\in[0,1],\\
u_-(-\infty,t)=x(t),u_+(+\infty,t)=y(t),\\
u_{\pm}(0,1)=\gamma(\pm\varepsilon)
\end{array}
\right.
\right\},
$$
(see figure below) and consider the moduli space
$$
\underline{{\mathcal M}}(x,y,H;f)=\big\{ (\varepsilon,u_{-},u_{+},\gamma)\,|\,(u_{-},u_{+},\gamma)\in{\mathcal M}_{\varepsilon}(x,y,H;f),\,\varepsilon\in[\varepsilon_0,\varepsilon_1]\big\},
$$
where $\varepsilon_0$ and $\varepsilon_1$ are fixed positive numbers.
\vskip5mm
\begin{center}\includegraphics[width=9cm,height=4.5cm]{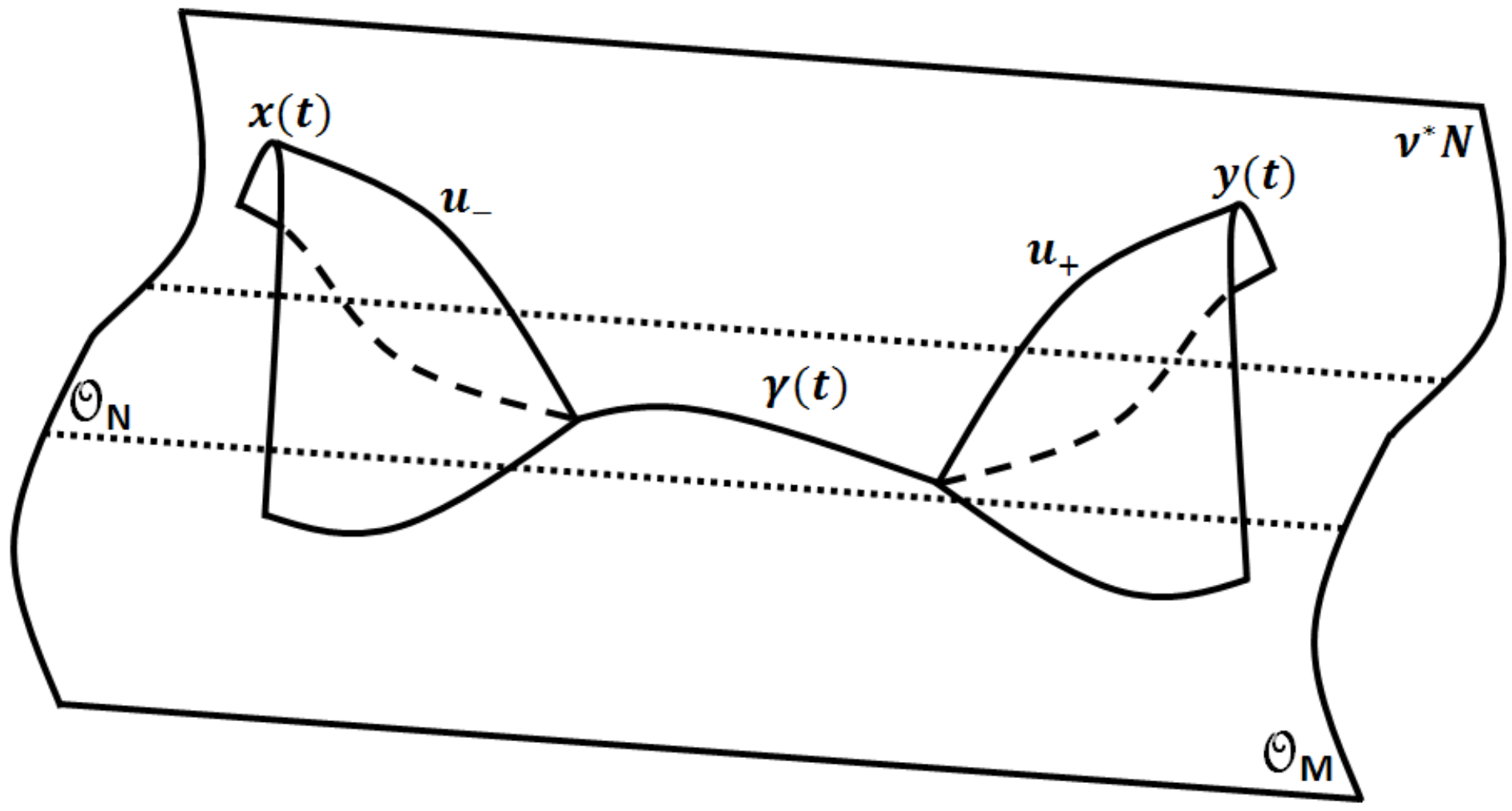}\\
Fig. 6. ${\mathcal M}_{\varepsilon}(x,y,H;f)$
\end{center}
\vskip5mm
For $\mu_N(y)=\mu_N(x)+1$, $\underline{{\mathcal M}}(x,y,H;f)$ is a zero--dimensional manifold. If $\mu_N(y)=\mu_N(x)$ then $\underline{{\mathcal M}}(x,y,H;f)$ is an one--dimensional manifold and we can describe its boundary.
\begin{prop}\label{boundary4}
Let $x,y\in CF_k(H)$. Then the topological boundary of $\underline{{\mathcal M}}(x,y,H;f)$ can be identified with
$$
\begin{aligned}
\partial \underline{{\mathcal M}}(x,y,H;f)=&{\mathcal M}_{\varepsilon_1}(x,y,H;f)\,\,\,\cup\,\,\,{\mathcal M}_{\varepsilon_0}(x,y,H;f)\\
&\cup\,\,\,\,\,\,\,\bigcup_{\mu_N(z)=\mu_N(x)-1}{\mathcal M}(x,z;H)\times\underline{\mathcal M}(z,y,H;f)\\
&\cup\,\,\,\,\,\,\,\bigcup_{\mu_N(z)=\mu_N(x)+1}\underline{\mathcal M}(x,z,H;f)\times{\mathcal M}(z,y;H).
\end{aligned}
$$
\end{prop}
\smallskip
\noindent {\it Proof:} Let us take a sequence $(\varepsilon_n,u_-^n,u_+^n,\gamma_n)\in\underline{{\mathcal M}}(x,y,H;f)$ that has no convergent subsequence in $W^{1,2}$--topology. Since a sequence $\varepsilon_n$ is bounded all uniform estimates for $u_\pm^n,\gamma_n$ hold uniformly on $\varepsilon$ (see Proposition \ref{boundary1}). Hence, sequences $u_-^n,u_+^n$ and $\gamma_n$ converge locally uniformly and $(u_-^n,u_+^n,\gamma_n)$ can break only once (for dimensional reason). The domain of $\gamma_n$ is bounded so trajectory $\gamma_n$ cannot break. The only remaining possibilities are:
\begin{itemize}
\item[(1)] There is a subsequence which converges to an element of ${\mathcal M}_{\varepsilon_1}(x,y,H;f)$ or ${\mathcal M}_{\varepsilon_0}(x,y,H;f)$.
\item[(2)] There is a subsequence which converges weakly to an element of ${\mathcal M}(x,z;H)\times\underline{\mathcal M}(z,y,H;f)$.
\item[(3)] There is a subsequence which converges weakly to an element of $\underline{\mathcal M}(x,z,H;f)\\
\times{\mathcal M}(z,y;H)$.\qed
\end{itemize}
Now, we define moduli space similar to $\overline{\mathcal M}(p,q,f;H)$, except that we are not using fixed Hamiltonian $H$ but a homotopy of Hamiltonians $H_{\delta}$, $0\le\delta\le1$, that connects given Hamiltonians $H_0$ and $H_1$,
$$
\overline{\mathcal M}(p,q,f;H_\delta)=\big\{(\delta,\gamma_-,\gamma_+,u)\,|\,(\gamma_-,\gamma_+,u)\in{\mathcal M}_{R_0}(p,q,f;H_\delta)),\,0\le\delta\le1\big\}.
$$
The dimension of this manifold is $m_f(p)-m_f(q)+1$ and its boundary is described in the following proposition.
\begin{prop}\label{boundary5}
Let $p,q\in CM_k(f)$. Then topological boundary of one--dimensional manifold $\overline{\mathcal M}(p,q,f;H_\delta)$ can be identified with
$$
\begin{aligned}
\partial\overline{\mathcal M}(p,q,f;H_\delta)=&{\mathcal M}_{R_0}(p,q,f;H_0)\,\,\,\cup\,\,\,{\mathcal M}_{R_0}(p,q,f;H_1)\\
&\cup\bigcup_{m_f(r)=k-1}{\mathcal M}(p,r;f)\times\overline{\mathcal M}(r,q,f;H_\delta)\\
&\cup\bigcup_{m_f(r)=k+1}\overline{\mathcal M}(p,r,f;H_\delta)\times{\mathcal M}(r,q;f).
\end{aligned}
$$
\end{prop}
\smallskip
\noindent {\it Proof:} Proof is essentially the same as for Proposition~\ref{boundary3}.\qed\\

So far, we have discussed moduli spaces defined by a family of Hamiltonians with a fixed Morse function $f$. It will be useful to consider moduli spaces similar to ${\mathcal M}(p,f;x,H)$, that depend on a family of Morse functions and a family of Hamiltonians. Let $(f^{\alpha\beta}_{s,\delta},H^{\alpha\beta}_{s,\delta})$, $0\leq\delta\leq1$, be a homotopy connecting $(f^\alpha,H^{\alpha\beta}_s)$ for $\delta=0$ and $(f^{\alpha\beta}_s,H^\beta)$ for $\delta=1$. Here, $f^{\alpha\beta}_s$ is a homotopy connecting two Morse functions $f^\alpha$ and $f^\beta$
$$
f^{\alpha\beta}_s=\begin{cases} f^\alpha, & s\le -T-1, \\ f^\beta, & s\ge -T. \end{cases}
$$
In the same way, $H^{\alpha\beta}_s$ is a homotopy connecting two Hamiltonians $H^\alpha$ and $H^\beta$
$$
H^{\alpha\beta}_s=\begin{cases} H^\alpha, & s\le T, \\ H^\beta, & s\ge T+1. \end{cases}
$$
We choose homotopy $(f^{\alpha\beta}_{s,\delta},H^{\alpha\beta}_{s,\delta})$ such that for any $\delta$ and $s$ negative (positive) enough, $f^{\alpha\beta}_{s,\delta}$ is equal to $f^\alpha$ ($H^{\alpha\beta}_{s,\delta}$ is equal to $H^\beta$).
Let
\begin{equation}\label{f,H_homot}
\widehat{{\mathcal M}}(p^{\alpha},f^{\alpha\beta}_{s,\delta};x^{\beta},H^{\alpha\beta}_{s,\delta})
=
\left\{ (\delta,\gamma,u) \left|
\begin{array}{ll}
\gamma:(-\infty,0]\rightarrow N, \\
u:[0,+\infty)\times[0,1]\rightarrow T^*M, \\
\frac{d\gamma}{ds}=-\nabla f^{\alpha\beta}_{s,\delta}(\gamma(s)), \\
\frac{\partial u}{\partial s}+J(\frac{\partial u}{\partial t}-X_{\rho_R^+H^{\alpha\beta}_{s,\delta}}(u))= 0, \\
\gamma(-\infty)=p^\alpha,\\
u(s,0)\in o_M, u(s,1)\in\nu^*N,s\ge 0,\\
u(0,t)\in o_M,t\in[0,1],\\
u(+\infty,t)=x^\beta(t),\\
\gamma(0)=u(0,1)
\end{array}
\right.
\right\}.
\end{equation}
The dimension of this manifold is $m_{f^{\alpha}}(p^\alpha)-(\mu_N(x^\beta)+\frac{1}{2}\dim N)+1$. The manifolds
$$
{\mathcal M}(p^{\alpha},f^{\alpha\beta}_s;x^{\beta},H^{\beta})
=
\left\{ (\gamma,u) \left|
\begin{array}{ll}
\gamma:(-\infty,0]\rightarrow N, \\
u:[0,+\infty)\times[0,1]\rightarrow T^*M, \\
\frac{d\gamma}{ds}=-\nabla f^{\alpha\beta}_s(\gamma(s)), \\
\frac{\partial u}{\partial s}+J(\frac{\partial u}{\partial t}-X_{\rho_R^+H^{\beta}}(u))= 0, \\
\gamma(-\infty)=p^\alpha,\\
u(s,0)\in o_M, u(s,1)\in\nu^*N,s\ge 0,\\
u(0,t)\in o_M,t\in[0,1],\\
u(+\infty,t)=x^\beta(t),\\
\gamma(0)=u(0,1)
\end{array}
\right.
\right\},
$$
and
\begin{equation}\label{H_homot}
{\mathcal M}(p^{\alpha},f^{\alpha};x^{\beta},H^{\alpha\beta}_s)
=
\left\{ (\gamma,u) \left|
\begin{array}{ll}
\gamma:(-\infty,0]\rightarrow N, \\
u:[0,+\infty)\times[0,1]\rightarrow T^*M, \\
\frac{d\gamma}{ds}=-\nabla f^{\alpha}(\gamma), \\
\frac{\partial u}{\partial s}+J(\frac{\partial u}{\partial t}-X_{\rho_R^+H^{\alpha\beta}_s}(u))= 0, \\
\gamma(-\infty)=p^\alpha,\\
u(s,0)\in o_M, u(s,1)\in\nu^*N,s\ge 0,\\
u(0,t)\in o_M,t\in[0,1],\\
u(+\infty,t)=x^\beta(t),\\
\gamma(0)=u(0,1)
\end{array}
\right.
\right\},
\end{equation}
are the two components of a boundary $\partial\widehat{{\mathcal M}}(p^{\alpha},f^{\alpha\beta}_{s,\delta};x^{\beta},H^{\alpha\beta}_{s,\delta})$ which we completely describe in the next proposition.

\begin{prop}\label{boundary6}
Let $m_{f^{\alpha}}(p^\alpha)=\mu_N(x^\beta)+\frac{1}{2}\dim N$. Then topological boundary of one--dimensional manifold $\widehat{{\mathcal M}}(p^{\alpha},f^{\alpha\beta}_{s,\delta};x^{\beta},H^{\alpha\beta}_{s,\delta})$ can be identified with
$$
\begin{aligned}
\partial\widehat{{\mathcal M}}(p^{\alpha},f^{\alpha\beta}_{s,\delta};x^{\beta},H^{\alpha\beta}_{s,\delta})=
&{\mathcal M}(p^{\alpha},f^{\alpha\beta}_s;x^{\beta},H^{\beta})
    \,\,\,\cup\,\,\,{\mathcal M}(p^{\alpha},f^{\alpha};x^{\beta},H^{\alpha\beta}_s)\\
&\cup\bigcup_{m_{f^{\alpha}}(q^{\alpha})=m_{f^{\alpha}}(p^{\alpha})-1}
    {\mathcal M}(p^{\alpha},q^{\alpha};f^{\alpha})\times
    \widehat{{\mathcal M}}(q^{\alpha},f_{s,\delta}^{\alpha\beta};x^{\beta},H_{s,\delta}^{\alpha\beta})\\
&\cup\,\,\,\bigcup_{\mu_N(y^\beta)=\mu_N(x^\beta)+1}
    \widehat{{\mathcal M}}(p^{\alpha},f_{s,\delta}^{\alpha\beta};y^{\beta},H_{s,\delta}^{\alpha\beta})
    \times
    {\mathcal M}(y^{\beta},x^{\beta};H^{\beta}).
\end{aligned}
$$
\end{prop}
\smallskip
\noindent {\it Proof:} Proof is essentially the same as for Proposition~\ref{boundary3}.\qed\\
\smallskip
\section{\textbf{Isomorphism}}\label{izomorfizam}
We saw in Proposition~\ref{boundary1} and Proposition~\ref{boundary2} that ${\mathcal M}(p,f;x,H)$ and ${\mathcal M}(x,H;p,f)$ are finite sets if $m_f(p)=\mu_N(x)+\frac{1}{2}\dim N$. Cardinal numbers of these sets (modulo 2) will be denoted by $n(p,f;x,H)$ and $n(x,H;p,f)$. Let us define homomorphisms on generators:
$$
\begin{aligned}
\phi:CF_k(H)\to CM_k(f),\,\,\,&\phi(x)=\sum_{m_f(p)=k}n(x,H;p,f)\,p,\\
\psi:CM_k(f)\to CF_k(H),\,\,\,&\psi(p)=\sum_{\mu_N(x)=k-\frac{1}{2}\dim N}n(p,f;x,H)\,x.
\end{aligned}
$$

\begin{prop}\label{chain_maps}
Homomorphisms $\phi$ and $\psi$ are well defined chain maps.
\end{prop}
\smallskip

\noindent {\it Proof:} It follows from Propositions~\ref{boundary1}, Propositions~\ref{boundary2} and from the way the chain complexes $CM_*(f)$ and $CF_*(H)$ are graded that these homomorphisms are well defined .\\
We prove that $(\phi\circ\partial_F-\partial_M\circ\phi)(x)=0$ for all $x\in CF_k(H)$.
$$\begin{aligned}
(\phi\circ\partial_F-\partial_M\circ\phi)(x)=&\sum_{m_f(q)=k-1}\bigg(\sum_{\mu_N(y)+\dim N/2=k-1}n(x,y;H)\,n(y,H;q,f)\bigg)q\,-\\
-&\sum_{m_f(q)=k-1}\bigg(\sum_{m_f(p)=k}n(x,H;p,f)\,n(p,q;f)\bigg)q.
\end{aligned}$$
Let $p\in CM_k(f)$, $q\in CM_{k-1}(f)$ and $y\in CF_{k-1}(H)$. From Proposition~\ref{boundary2} it follows
$$
\sum_{\mu_N(y)+\dim N/2=k-1}n(x,y;H)\,n(y,H;q,f)-\sum_{m_f(p)=k}n(x,H;p,f)\,n(p,q;f)=0,
$$
since it is (modulo 2) number of ends of one--dimensional manifold ${\mathcal M}(x,H;q,f)$. So, $(\phi\circ\partial_F-\partial_M\circ\phi)(x)=0$. The proof of identity $\psi\circ\partial_M=\partial_F\circ\psi$ is analogous.
\qed
\\

From the previous proposition it follows that $\phi$ and $\psi$ induce homomorphisms in homology, $$\Phi:HF_k(H)\rightarrow HM_k(f),\,\,\,\Psi:HM_k(f)\rightarrow HF_k(H).$$
These homomorphisms are PSS--type isomorphisms. Now, we can prove Theorem~\ref{main_thm}. From the fact that these homomorphisms are inverse to each other it will immediately follow that $\Phi$ and $\Psi$ are isomorphisms. In order to show that we prove that $\phi\circ\psi$ and $\psi\circ\phi$ are maps chain homotopic to the identity.\\

\noindent {\it Proof of Theorem~\ref{main_thm}:} If we look at a composition of homomorphisms $\phi$ and $\psi$, $$\phi\circ\psi(p)=\sum_{m_f(q)=k}\bigg(\sum_{\mu_N(x)+\dim N/2=k}n(p,f;x,H)\,n(x,H;q,f)\bigg)q,$$ we can see that $\sum_{x}n(p,f;x,H)n(x,H;q,f)$ is number of points of a set $$\cup_{x}{\mathcal M}(p,f;x,H)\times{\mathcal M}(x,H;q,f),$$ which is a component of boundary $\partial\overline{\mathcal M}(p,q,f;H)$.\\
Similarly to~\cite{KM} we define homomorphisms $l$ and $j$,
$$
\begin{aligned}
l:CM_k(f)\rightarrow CM_k(f),\,\,\,\,\,\,\,\,&l(p)=\sum_{m_f(q)=k}n(p,q,f;H)\,q,\\
j:CM_k(f)\rightarrow CM_{k+1}(f),\,\,\,&j(p)=\sum_{m_f(r)=k+1}\overline{n}(p,r,f;H)\,r.
\end{aligned}
$$
Here $n(p,q,f;H)$ is the number of intersections of a space of perturbed holomorphic discs with the unstable manifold $W^u(p,f)$ and the stable manifold $W^s(q,f)$. We consider discs with half of a boundary on the zero section, $o_M$, and half of a boundary on the conormal bundle, $\nu^*N$. In other words, $n(p,q,f;H)$ is the number of elements of ${\mathcal M}_{R_0}(p,q,f;H)$. By $\overline{n}(p,r,f;H)$ we denote the number of elements of a zero--dimensional manifold $\overline{\mathcal M}(p,r,f;H)$. A sum
$$
\sum_{m_f(r)=k-1}n(p,r;f)\,\overline{n}(r,q,f;H)
$$
corresponds to a sum that occurs in $j\circ\partial_M$, and
$$
\sum_{m_f(r)=k+1}\overline{n}(p,r,f;H)\,n(r,q;f)
$$
corresponds to a sum in $\partial_M\circ j$. From Proposition~\ref{boundary3} follows
$$
\phi\circ\psi-l=\partial_M\circ j+j\circ\partial_M.
$$
Now, we prove that homomorphism in homology,
$$
L:HM_k(f)\rightarrow HM_k(f),
$$
induced by chain homomorphism $l$ does not depend on Hamiltonian $H$. Let $H_0$ and $H_1$ be two Hamiltonians and $H_\delta$, $0\le \delta\le 1$, a homotopy between them. $l_0$ and $l_1$ are chain homomorphisms corresponding to the $H_0$ and $H_1$. From Proposition~\ref{boundary5} we get the relation
$$
l_1-l_0=\partial_M\circ j_\delta+j_\delta\circ\partial_M,
$$
where
$$
j_\delta:CM_k(f)\rightarrow CM_{k+1}(f),\hskip5mm j_\delta(p)=\sum_{m_f(r)=k+1}\overline{n}(p,r,f;H_\delta)\,r.
$$
Here, $\overline{n}(p,r,f;H_\delta)$ is the number of elements of $\overline{\mathcal M}(p,r,f;H_\delta)$.
If we choose homotopy between our Hamiltonian $H$ and 0 we conclude that a map $l$ is chain homotopic to a map $i:CM_k(f)\to CM_k(f)$, $$
i(p)=\sum_{m_f(q)=k}n(p,q,f;0)q.
$$
Thus, $L$ and a map $I$, induced by $i$, are the same maps in homology. We explained above that unperturbed holomorphic disc with half of a boundary on the zero section and the other half on the conormal bundle is constant. It follows that $n(p,q,f;0)$ is the number of points in $W^u(p,f)\cap W^s(q,f)$. Considering Morse indices of $p$ and $q$ we get $I=\mathbb{Id}$.

We use the same idea to prove $\Psi\circ\Phi=\mathbb {Id}$. The composition $\psi\circ\phi$ is chain homotopic to some chain homomorphism $r:CF_k(H)\to CF_k(H)$ which induces the identity in homology. If we denote by $n_\varepsilon(x,y,H;f)$ the number of elements of a zero--dimensional manifold ${\mathcal M}_{\varepsilon}(x,y,H;f)$ then the map analogous to $l$ is
$$
r(x)=\sum_{\mu_N(y)=\mu_N(x)}n_\varepsilon(x,y,H;f)\,y.
$$
Similarly to the first part of the proof, a homomorphism in homology induced by $r$ is independent of the choice of $\varepsilon$. Let $r_0$ and $r_1$ be homomorphisms corresponding to the values $\varepsilon_0$ and $\varepsilon_1$. We define a chain homomorphism
$$
s:CF_k(H)\rightarrow CF_{k+1}(H),\hskip5mm s(x)=\sum_{\mu_N(y)+\dim N/2=k+1}\underline{n}(x,y,H;f)\,y,
$$
where $\underline{n}(x,y,H;f)$ denotes the number of elements of $\underline{{\mathcal M}}(x,y,H;f)$. From Proposition~\ref{boundary4}
we conclude
$$
r_0-r_1=s\circ\partial_F+\partial_F\circ s.
$$
If we pass to the limit as $\varepsilon\to0$ we get that $\psi\circ\phi$ if chain homotopic to a homomorphism $\widetilde{i}:CF_k(H)\rightarrow CF_k(H)$,
$$
\widetilde{i}(x)=\sum_{\mu_N(y)+\dim N/2=k}\widetilde{n}(x,y;H)\,y,
$$
where $\widetilde{n}(x,y;H)$ is the number of elements of a zero--dimensional manifold
$$
\widetilde{\mathcal M}(x,y;H)=\left\{ (u_{-},u_{+}) \left|
\begin{array}{ll}
u_{-}:(-\infty,0]\times[0,1]\rightarrow T^*M,\\
u_{+}:[0,+\infty)\times[0,1]\rightarrow T^*M, \\
\frac{\partial u_{\pm}}{\partial s}+J(\frac{\partial u_{\pm}}{\partial t}-X_{\rho_R^{\pm}H}(u_{\pm}))= 0, \\
u_{\pm}(\pm s,0)\in o_M, u_{\pm}(\pm s,1)\in\nu^*N,s\ge0,\\
u_{\pm}(0,t)\in o_M, t\in[0,1],\\
u_-(-\infty,t)=x(t),u_+(+\infty,t)=y(t),\\
u_{+}(0,1)=u_{-}(0,1)
\end{array}
\right.
\right\},
$$
(see figure below). The rest of this proof is dedicated to showing that counting the number of elements of $\widetilde{\mathcal M}(x,y;H)$ is the same as counting the pseudo holomorphic strips between $x$ and $y$ (at the homology level). The main idea is to show that $\widetilde{\mathcal M}$ is cobordant to the manifold that consists of appropriate pseudo holomorphic strips.
\vskip5mm
\begin{center}\includegraphics[width=9cm,height=4.5cm]{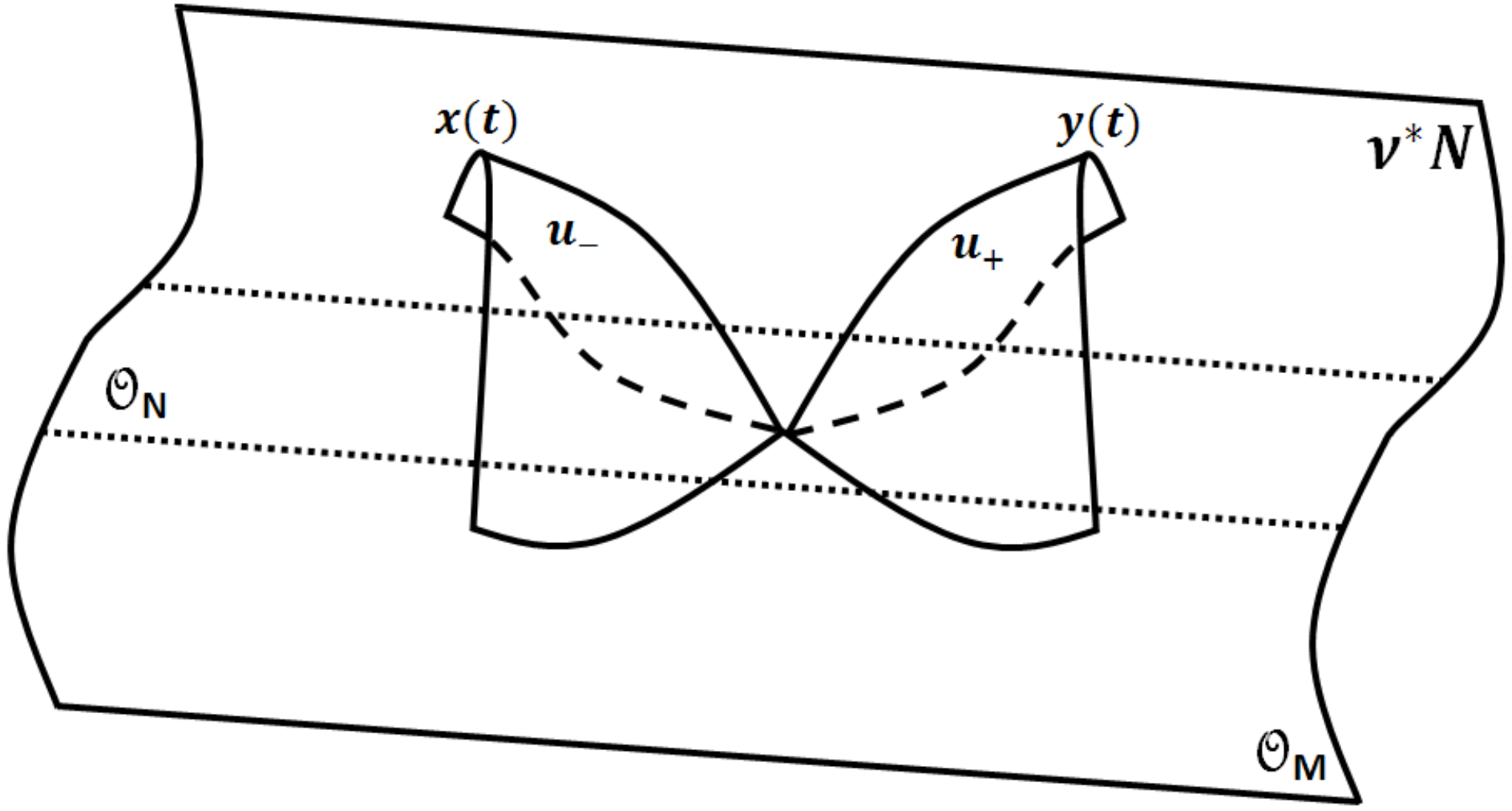}\\
Fig. 7. $\widetilde{\mathcal M}(x,y;H)$
\end{center}
\vskip5mm

It is more convenient to deal with the whole strip as a domain of $u_\pm$ instead the half--strip (for example, when we preform gluing we have a translation along $s$--axis). Instead of $\widetilde{\mathcal M}(x,y;H)$ we can observe the manifold (we keep the same notation)
$$
\widetilde{\mathcal M}(x,y;H)=\left\{ (u_{-},u_{+}) \left|
\begin{array}{ll}
u_{\pm}:\mathbb{R}\times[0,1]\rightarrow T^*M,\\
E(u_\pm)<\infty,\\
\frac{\partial u_{\pm}}{\partial s}+J(\frac{\partial u_{\pm}}{\partial t}-X_{\rho_R^{\pm}H}(u_{\pm}))= 0, \\
u_{\pm}(s,0)\in o_M, u_{\pm}(s,1)\in\nu^*N,s\in\mathbb{R},\\
u_-(-\infty,t)=x(t),u_+(+\infty,t)=y(t),\\
u_{-}(+\infty)=u_{+}(-\infty)
\end{array}
\right.
\right\}.
$$
Smooth functions $\rho^\pm_R$ are now defined on the whole $\mathbb{R}$ as follows
$$
\rho_R^+(s)=\begin{cases} 1, & s\ge R+1, \\ 0, & s\le R, \end{cases}
$$
and
$$
\rho_R^-(s)=\begin{cases} 1, & s\le -R-1, \\ 0, & s\ge -R. \end{cases}
$$
If we use the second definition of $\widetilde{\mathcal M}$, strip $u_-$ ($u_+$) is holomorphic for $s\ge-R$ ($s\le R$) and has finite energy. So, $u_\pm$ admits a unique continuous extension $u_\pm(\mp\infty)$ (see Section 4.5 in~\cite{MS1} and Theorem 3.1 in~\cite{FelS}). These extensions are points that belong to $N\subset o_M$ and we can omit the second argument of $u_\pm(\mp\infty)$. We can transform the domain of $u_\pm$ because the half--strip is conformally equivalent to a disc without a boundary point. From Remark~\ref{napomena} we know that the strip $\mathbb{R}\times[0,1]$ is conformal to $\mathbb{D}^2\setminus\{\pm1\}$. When we remove the singularity at one infinity point ($+\infty$ of $u_-$ and $-\infty$ of $u_+$) we obtain a disc without one boundary point.

We define auxiliary manifolds
 $${\mathcal M}_R(x,y;H)=\left\{ u \left|
\begin{array}{ll}
u:\mathbb{R}\times[0,1]\rightarrow T^*M,\\
\frac{\partial u}{\partial s}+J(\frac{\partial u}{\partial t}-X_{\rho_RH}(u))= 0, \\
u(s,0)\in o_M, u(s,1)\in\nu^*N,s\in{\mathbb R},\\
u(-\infty,t)=x(t),u(+\infty,t)=y(t)
\end{array}
\right.
\right\}
$$
and
$$\check{{\mathcal M}}(x,y;H)=\left\{ (R,u) \left| R\ge R_0,\,u\in{\mathcal M}_R(x,y;H)
\right.
\right\}.
$$
Here, $\rho_R:{\mathbb R}\rightarrow{\mathbb R}$ is a smooth function such that
$$
\rho_R(s)=\begin{cases} 1, & |s|\ge R+1, \\ 0, & |s|\le R. \end{cases}
$$
The boundary of the manifold $\check{{\mathcal M}}(x,y;H)$ can be identified with
\begin{equation}\label{theLast1}
\begin{aligned}
\partial\check{{\mathcal M}}={\mathcal M}_{R_0}(x,y;H)\,\,&\cup\,\,\widetilde{\mathcal M}(x,y;H)\\
&\cup\,\,\bigcup_{z}\mathcal{M}(x,z;H)\times\check{{\mathcal M}}(z,y;H)\\
&\cup\,\,\bigcup_{z}\check{{\mathcal M}}(x,z;H)\times{\mathcal M}(z,y;H).
\end{aligned}
\end{equation}
Now we explain the equality~(\ref{theLast1}). It is clear how the last two terms on the right--hand side appear at the boundary. Elements from ${\mathcal M}_{R_0}(x,y;H)$ appear at the boundary when $R_n\to R_{0}$. The most complicated part is to prove that $\widetilde{\mathcal M}$ is a part of the boundary. We show that in two steps. In Step A we explain why it holds $\partial\check{{\mathcal M}}\subset\widetilde{\mathcal M}$. And in Step B we show the opposite inclusion, $\widetilde{\mathcal M}\subset\partial\check{{\mathcal M}}$.\\

{\it Step A.} When $R_n\to+\infty$ we can identify the limit of $u_n\in{\mathcal M}_{R_n}(x,y;H)$ with the element from $\widetilde{\mathcal M}(x,y;H)$ using the reparametrization:
$$
u_n^-(s,t)=u_n(s-R_n+R_0,t),\,\,\,u_n^+(s,t)=u_n(s+R_n-R_0,t).
$$
Strip $u_n^-$ satisfies the equation
$$
\frac{\partial u_n^-}{\partial s}+J\left(\frac{\partial u_n^-}{\partial t}-X_{\rho_{R_n-}H}(u_n^-)\right)=0,
$$
and the boundary conditions
$$
\begin{aligned}
  &u_n^-(s,0)=u_n(s-R_n+R_0,0)\in o_M, \\
  &u_n^-(s,1)=u_n(s-R_n+R_0,1)\in\nu^*N,
\end{aligned}
$$
for $s\in{\mathbb R}$.
Function $\rho_{R_n-}$ is equal to
$$
\rho_{R_n-}(s)=\begin{cases} 0, &-R_0\le s\le 2R_n-R_0, \\ 1, & s\in (-\infty,-R_0-1]\cup[2R_n-R_0+1,+\infty). \end{cases}
$$
Positive strip $u_n^+$ also satisfies perturbed Cauchy–-Riemann equation
$$
\frac{\partial u_n^+}{\partial s}+J\left(\frac{\partial u_n^+}{\partial t}-X_{\rho_{R_n+}H}(u_n^+)\right)=0,
$$
the line $u_n^+({\mathbb R}\times\{0\})$ is on the zero section and $u_n^+({\mathbb R}\times\{1\})$ is on the conormal bundle. Function $\rho_{R_n+}$ is defined as
$$
\rho_{R_n+}(s)=\begin{cases} 0, &-2R_n+R_0\le s\le R_0, \\ 1, & s\in (-\infty,-2R_n+R_0-1]\cup[R_0+1,+\infty). \end{cases}
$$
Strip $u_n^\pm$ converges locally uniformly with all derivatives to some $u^\pm$ that satisfies the equation
$$
\frac{\partial u^\pm}{\partial s}+J\left(\frac{\partial u^\pm}{\partial t}-X_{\rho^\pm_{R_0} H}(u^\pm)\right)=0.
$$
It is obvious that $u^-(-\infty,t)=x(t)$ and $u^+(+\infty,t)=y(t)$. At the $+\infty$--end the strip $u^-$ converges to some point $p\in N\subset o_M$ since $u^-$ is holomorphic for $s\ge -R_0$ and it has finite energy. The positive strip $u^+$ is holomorphic at the $-\infty$ and it converges to some point $q\in N$. Since
$$
u_n^-(R_n-R_0,t)=u_n^+(-R_n+R_0,t)
$$
we conclude that $p=q$. Thus, the pair $(u^-,u^+)$ belongs to $\widetilde{\mathcal M}(x,y;H)$.\\

{\it Step B.} For a given $(u_-,u_+)\in\widetilde{\mathcal M}$ we can find a sequence of elements $(R,\omega_R)\in\check{{\mathcal M}}$ that Gromov converges to $(u_-,u_+)$ as $R\to+\infty$ (see Theorem 4.1.2 in~\cite{BC}, Chapter 4.7 in~\cite{Urs} and Theorem 7.1 in~\cite{FelS}).

The main technique is gluing and goes as follows. Strips $u_-$ and $u_+$ are holomorphic around point $u_-(+\infty)=u_+(-\infty)$ and we can preglue them to obtain a map $u_R$. This is an approximate solution of the Cauchy–-Riemann, $u_R$ satisfies this equation everywhere besides the small neighbourhood of $u_R(0)=u_-(+\infty)=u_+(-\infty)$. After this we construct a right inverse to the linearization $D_{u_R}$ of the operator $\overline{\partial}$. Using an implicit function theorem we find a genuine solution to this equation, $\omega_R$, that is in a neighborhood of an approximate solution.

Biran and Cornea in~\cite{BC} glued two holomorphic discs with the boundary on one Lagrangian submanifold. Frauenfelder in~\cite{Urs} and Schm\"{a}schke in~\cite{FelS} worked with two cleanly intersecting (compact) submanifolds in a compact symplectic manifold. The cotangent bundle is not a compact manifold but, with appropriate choice on an almost complex structure (see the definition of $j^c$ below), the image of every holomorphic strip lies in a compact subset of $T^*M$ (see Theorem 3.2. in~\cite{O2}). So we can assume that everything happens in the compact subset of our symplectic manifold. We also need special Riemannian metrics on $T^*M$ such that $o_M$ and $\nu^*N$ are totally geodesic submanifolds with respect to these metrics.

Following the~\cite{O2},~\cite{Urs1} and~\cite{Sim} we explain choices on almost complex structures and Riemannian metrics on $T^*M$. Fix a Riemannian metric $g$ on $M$. The associated Levi-Civita connection induces the canonical almost complex structure on $T^*M$, which we denote by $J_g$. We define the subset $j^c$ of the set of almost complex structure on $T^*M$
$$j^c=\{J\,|\,J\text{ is compatible to }\omega,\,J=J_g\text{ outside a compact subset in } T^*M\}.$$
Let $J_t$ be a smooth path in $j^c$. Then there exists a smooth family of metrics $g_t$ such that
\begin{enumerate}
  \item $o_M$ is totally geodesic with respect to $g_0$ and $J_0(q)T_qo_M$ is the orthogonal complement of $T_qo_M$ for every $q\in o_M$,
  \item $\nu^*N$ is totally geodesic with respect to $g_1$ and $J_1(q)T_q(\nu^*N)$ is the orthogonal complement of $T_q(\nu^*N)$ near the intersection point of two holomorphic strips that we glue,
  \item $g_t(J_t(q)u,J_t(q)v)=g_t(u,v)$ for $q\in T^*M$ and $u,v\in T_q(T^*M)$.
\end{enumerate}
We can define a metric $g_0$ such that
\begin{enumerate}
  \item $o_M$ is totally geodesic with respect to $g_0$ and $J_0(q)T_qo_M$ is the orthogonal complement of $T_qo_M$ for every $q\in o_M$,
  \item $g_0(J_0(q)u,J_0(q)v)=g_0(u,v)$ for $q\in T^*M$ and $u,v\in T_q(T^*M)$,
\end{enumerate}
see~\cite{Urs1} for details. In the same way we can define a metric $g_1$ that satisfies the same properties for the submanifold $\nu^*N$. In~\cite{Urs1} the author assumes that the Lagrangian submanifold is compact. The conormal bundle $\nu^*N$ is not a compact manifold, in general. But it is enough to find a metric such that $\nu^*N$ is a totally geodesic submanifold near the $N\subset o_M$, not the whole $\nu^*N$. The linear combination
$$g_t(u,v)=\overline{g}_t(u,v)+\overline{g}_t(J_tu,J_tv),
$$
where $\overline{g}_t(u,v)=(1-t)g_0(u,v)+tg_1(u,v)$, gives us the appropriate family of metrics (see also~\cite{Sim}).

All the other technical details of gluing are the same as in~\cite{FelS}.
\\

Now, we return to the homomorphism $\widetilde{i}$. Using the one--dimensional component of $\check{{\mathcal M}}(x,y;H)$ and the description of its boundary~(\ref{theLast1}) we conclude that $\widetilde{i}$ (i.e. $\psi\circ\phi$) is chain homotopic to a map $$k:x\mapsto\sum_{\mu_N(y)=\mu_N(x)}n(x,y;H)\,y.$$ If there is a non--constant holomorphic strip that connects Hamiltonian orbits $x$ and $y$ then $\mu_N(x)>\mu_N(y)$.  It follows that $$
n(x,y;H)=\begin{cases} 1, & x=y, \\ 0, & x\neq y, \end{cases}
$$
i.e. the map $k$ induces the identity in homology $HF_*(H)$.\qed

\smallskip

\section{\textbf{Commutative diagram}}\label{dijagram}

\noindent {\it Proof of Theorem~\ref{diagram_commutes}:} This theorem states that
$$S^{\alpha\beta}\circ\Psi^{\alpha}=\Psi^\beta\circ T^{\alpha\beta}.
$$
Composition on the left--hand side is generated with map $\sigma^{\alpha\beta}\circ\psi^\alpha$, and the right--hand side is generated with $\psi^\beta\circ\tau^{\alpha\beta}$ on a chain level. Idea is to prove that these maps on chain level are homotopic to each other.

We separate proof in two steps. In Step 1 and Step 2 we define new maps $\chi$ and $\xi$ that are homotopic to $\sigma^{\alpha\beta}\circ\psi^\alpha$ and $\psi^\beta\circ\tau^{\alpha\beta}$, respectively. In the conclusion of the proof we show that $\chi$ and $\xi$ are chain homotopic maps.

{\it Step 1.} From definitions it follows
$$
(\sigma^{\alpha\beta}\circ\psi^\alpha)(p^\alpha)=\sum_{x^\alpha,x^\beta}n(p^\alpha,f^\alpha;x^\alpha,H^\alpha)\,n(x^\alpha,x^\beta;H^{\alpha\beta})\,x^\beta.
$$
It means that $\sigma^{\alpha\beta}\circ\psi^\alpha$ counts the number of points of a set
$$
\bigcup_{x^\alpha}{\mathcal M}(p^\alpha,f^\alpha;x^\alpha,H^\alpha)\times{\mathcal M}(x^\alpha,x^\beta;H^{\alpha\beta}),
$$
where ${\mathcal M}(x^\alpha,x^\beta;H^{\alpha\beta})$ denotes the set of solutions of~(\ref{eq:2}). Summation is taken over $x^\alpha,\,x^\beta$ such that
$$m_{f^\alpha}(p^\alpha)=\mu_N(x^\alpha)+\frac{1}{2}\dim N=\mu_N(x^\beta)+\frac{1}{2}\dim N.
$$
Let us define a family of homotopies between Hamiltonians $H^\alpha$ and $H^\beta$:
$$
H^{\alpha\beta}_{T,s}=\begin{cases} H^\alpha, & s\le T, \\ H^\beta, & s\ge T+1. \end{cases}
$$
We consider moduli space
$$
\breve{{\mathcal M}}(p^{\alpha},f^\alpha;x^{\beta},H^{\alpha\beta}_{T,s})
=
\left\{ (T,\gamma,u) \left|
\begin{array}{ll}
T\ge T_0,\\
\gamma:(-\infty,0]\rightarrow N, \\
u:[0,+\infty)\times[0,1]\rightarrow T^*M, \\
\frac{d\gamma}{ds}=-\nabla f^\alpha(\gamma), \\
\frac{\partial u}{\partial s}+J(\frac{\partial u}{\partial t}-X_{\rho_R^+H^{\alpha\beta}_{T,s}}(u))= 0, \\
\gamma(-\infty)=p^\alpha,\\
u(s,0)\in o_M, u(s,1)\in\nu^*N,s\ge 0,\\
u(0,t)\in o_M,t\in[0,1],\\
u(+\infty,t)=x^\beta(t),\\
\gamma(0)=u(0,1)
\end{array}
\right.
\right\}.
$$
Using the same idea as in the proof of Theorem~\ref{main_thm}, from gluing and compactness arguments it follows that boundary of $\breve{{\mathcal M}}$ can be described as
$$
\begin{aligned}
\partial\breve{{\mathcal M}}(p^{\alpha},f^\alpha;x^{\beta},H^{\alpha\beta}_{T,s})=
&{\mathcal M}(p^{\alpha},f^\alpha;x^{\beta},H^{\alpha\beta}_{T_0,s})\\
&\cup\bigcup_{x^\alpha}{\mathcal M}(p^\alpha,f^\alpha;x^\alpha,H^\alpha)\times
    {\mathcal M}(x^\alpha,x^\beta;H^{\alpha\beta})\\
&\cup\bigcup_{q^\alpha}
    {\mathcal M}(p^{\alpha},q^{\alpha};f^{\alpha})\times
    \breve{{\mathcal M}}(q^{\alpha},f^\alpha;x^{\beta},H_{T,s}^{\alpha\beta})\\
&\cup\bigcup_{y^\beta}
    \breve{{\mathcal M}}(p^{\alpha},f^\alpha;y^{\beta},H_{T,s}^{\alpha\beta})
    \times
    {\mathcal M}(y^{\beta},x^{\beta};H^{\beta}).
\end{aligned}
$$
The first element in a previous union is already described in~(\ref{H_homot}) (for fixed homotopy $H^{\alpha\beta}_s=H^{\alpha\beta}_{T_0,s}$). We define a map $\chi$ that counts the number of elements in ${\mathcal M}(p^{\alpha},f^\alpha;x^{\beta},H^{\alpha\beta}_s)$
$$\chi(p^\alpha)=\sum_{x^\beta}n(p^{\alpha},f^\alpha;x^{\beta},H^{\alpha\beta}_s)\,x^\beta.
$$
From description of topological boundary of $\breve{{\mathcal M}}$ we conclude that $\chi$ and $\sigma^{\alpha\beta}\circ\psi^\alpha$ are chain homotopic maps.

{\it Step 2.} Other composition satisfies the equation
$$\psi^\beta\circ\tau^{\alpha\beta}(p^\alpha)=\sum_{p^\beta,x^\beta}n(p^\alpha,p^\beta;f^{\alpha\beta})\,n(p^\beta,f^\beta;x^\beta,H^\beta)\,x^\beta.
$$
Now, $\psi^\beta\circ\tau^{\alpha\beta}$ counts the number of points of a set
$$
\bigcup_{p^\beta}{\mathcal M}(p^\alpha,p^\beta;f^{\alpha\beta})\times{\mathcal M}(p^\beta,f^\beta;x^\beta,H^\beta),
$$
where ${\mathcal M}(p^\alpha,p^\beta;f^{\alpha\beta})$ is the set of solutions of~(\ref{eq:parmorse}). Here, we take a sum over $p^\beta,\,x^\beta$ such that
$$m_{f^\alpha}(p^\alpha)=m_{f^\beta}(p^\beta)=\mu_N(x^\beta)+\frac{1}{2}\dim N.
$$
We define a map $\xi$ that counts the number of points in ${\mathcal M}(p^\alpha,f^{\alpha\beta}_s;x^\beta,H^\beta)$. It follows, similarly as in Step 1, that $\xi$ and $\psi^\beta\circ\tau^{\alpha\beta}$ are chain homotopic maps.

Using the moduli space $\widehat{{\mathcal M}}(p^{\alpha},f^{\alpha\beta}_{s,\delta};x^{\beta},H^{\alpha\beta}_{s,\delta})$, that we define in~(\ref{f,H_homot}), we prove $\chi$ and $\xi$ are chain homotopic maps. Let us define a chain homomorphism
$$
\begin{aligned}
&j:CM_{k-1}(f^\alpha)\to CF_k(H^\beta),\\
&j(p^\alpha)=\sum_{\mu_N(x^\beta)+\dim N/2=k}\widehat{n}(p^{\alpha},f^{\alpha\beta}_{s,\delta};x^{\beta},H^{\alpha\beta}_{s,\delta})\,x^\beta,
\end{aligned}
$$
where $\widehat{n}(p^{\alpha},f^{\alpha\beta}_{s,\delta};x^{\beta},H^{\alpha\beta}_{s,\delta})$ is the number of elements of a zero--dimensional manifold $\widehat{{\mathcal M}}(p^{\alpha},f^{\alpha\beta}_{s,\delta};x^{\beta},H^{\alpha\beta}_{s,\delta})$. From Proposition~\ref{boundary6} it follows that $$\xi-\chi+j\circ\partial_M+\partial_F\circ j=0.$$

\qed

\section{\textbf{Product in homology}}\label{products}
In this section we define a product
$$
\ast:HF_*(o_M,\nu^*N:H_1)\otimes HF_*(o_M,\nu^*N:H_2)\to HF_*(o_M,\nu^*N:H_3),
$$
and prove the subadditivity of spectral invariants with respect to this product.

We define a Riemannian surface with boundary $\Sigma$ as a disjoint union
$$
{\mathbb R}\times[-1,0]\sqcup{\mathbb R}\times[0,1]
$$
with identification $(s,0^-)\sim(s,0^+)$ for $s\geq0$ (see figure below). The surface $\Sigma$ is conformally equivalent to a closed disc with three boundary punctures. The complex structure on $\Sigma\backslash\{(0,0)\}$ is induced by the inclusion in ${\mathbb C}$, $(s,t)\mapsto s+it$. The complex structure at the point $(0,0)$ is given by the square root.
\vskip5mm
\begin{center}\includegraphics[width=6cm,height=2cm]{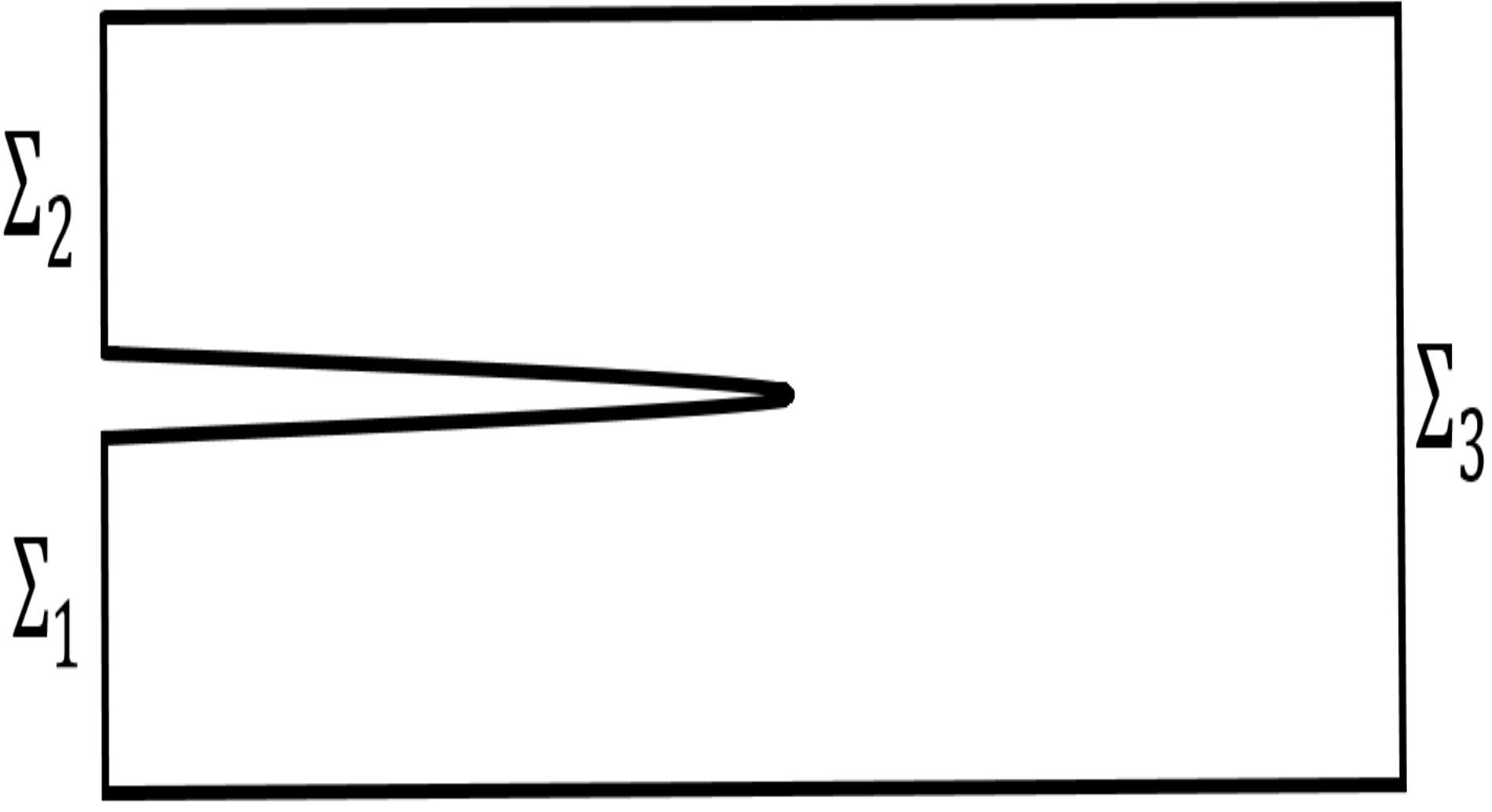}\\
Fig. 8. \textit{Riemannian surface $\Sigma$}
\end{center}
\vskip5mm

Denote by $\Sigma_1^-$, $\Sigma_2^-$, $\Sigma^+$ the two "incoming" and one "outgoing" ends, such that
 $$\begin{aligned}
 &\Sigma_1^-,\Sigma_2^-\approx [0,1]\times(-\infty,0],\\
 &\Sigma^+\approx [0,1]\times[0,+\infty).\\
 \end{aligned}$$
By $u_j^-:=u|_{\Sigma_j^-}$, $j=1,2$, and $u^+:=u|_{\Sigma^+}$ we denote a restriction of a map defined on the surface $\Sigma$.
Let $\rho^\pm:{\mathbb R}\to[0,1]$ denote the smooth cut--off functions such that
$$\rho^-(s)=\begin{cases}1,&s\le -2,\\0,&s\ge -1\end{cases}\quad\rho^+(s):=\rho^-(-s).$$
For $x_1^-\in CF_*(o_M,\nu^*N:H_1)$, $x_2^-\in CF_*(o_M,\nu^*N:H_2)$ and $x^+\in CF_*(o_M,\nu^*N:H_3)$ we define the moduli space
$$\begin{aligned}&{\mathcal M}(x_1^-,x_2^-;x^+)=&\left\{u:\Sigma\to T^*M\left|\begin{array}{l}
 \partial_su_j^-+J(\partial_tu_j^--X_{\rho^-H_j}\circ u_j^-)=0,\,j=1,2,\\
 \partial_su^++J(\partial_tu^+-X_{\rho^+H_3}\circ u^+)=0,\\
 \partial_su+J\partial_tu=0,\,\mbox{on}\;\Sigma_0:=\Sigma\setminus(\Sigma_1\cup \Sigma_2\cup \Sigma_3),\\
 u(s,-1)\in o_M,\,u(s,1)\in\nu^*N,\,s\in{\mathbb R},\\
 u(s,0^-)\in\nu^*N,\,u(s,0^+)\in o_M,\,s\le 0,\\
 u_j^-(-\infty,t)=x_j^-(t),\,j=1,2,\\
 u^+(+\infty,t)=x^+(t).
 \end{array}\right.\right\}.\end{aligned}$$

We use the notation
\begin{equation}\label{eq:perturbe}
\overline{\partial}_{J,H}(u)=0
\end{equation}
for perturbed Cauchy--Riemann equation that we consider in ${\mathcal M}(x_1^-,x_2^-;x^+)$. Elements of a moduli space ${\mathcal M}(x_1^-,x_2^-;x^+)$ are perturbed holomorphic discs $u$. The boundary of $u$ is on a Lagrangian submanifold $o_M\cup\nu^*N$ with clean self--intersection along $N$. This disc $u$ has one switch on a slit of the pants in the sense of~\cite{CTL} (or a jump in the sense of~\cite{AS}).

For generic choices of Hamiltonians and an almost complex structure, ${\mathcal M}(x_1^-,x_2^-;x^+)$ is a smooth manifold of finite dimension.

We give some more details on Fredholm analysis for this case. Let us define suitable Banach frame for Fredholm analysis. For $p>2$ we define
\begin{equation}
\mathcal{P}^{1,p}(x_1^-,x_2^-;x^+)=\left\{
\begin{array}{ll}
u\in W^{1,p}_{loc}(\Sigma,T^*M), \\
(\exists T>0)(\exists \xi_i^-\in W^{1,p}((-\infty,-T]\times[0,1],(x_i^-)^*T(T^*M))\\
\hskip10mm (\exists \xi^+\in W^{1,p}([T,+\infty)\times[0,1],(x^+)^*T(T^*M))\\
\hskip10mm u_i^-(s,t)=exp_{x_i^-(t)}\xi_i^-(s,t),s\le -T,\,i\in\{1,2\},\\
\hskip10mm u^+(s,t)=exp_{x^+(t)}\xi^+(s,t),s\ge T,\\
u(s,-1)\in o_M, u(s,1)\in\nu^*N,\,s\in{\mathbb{R}},\\
u(s,0^-)\in \nu^*N, u(s,0^+)\in o_M,\, s\leq0,\\
\lim\limits_{s\to-\infty}u_i^-(s,t)=x_i^-(t),\,i\in\{1,2\},\\
\lim\limits_{s\to+\infty}u^+(s,t)=x^+(t).
\end{array}
\right.
\end{equation}
$\mathcal{P}^{1,p}(x_1^-,x_2^-;x^+)$ is a Banach manifold and charts are obtained using exponential map. For $u\in\mathcal{P}^{1,p}$ it holds
$$
T_u\mathcal{P}^{1,p}(x_1^-,x_2^-;x^+)=W^{1,p}_\Lambda(u^*T(T^*M)),
$$
where on the right--hand side we have $W^{1,p}$-section of a vector bundle $u^*T(T^*M)\to\Sigma$ with Lagrangian boundary conditions:
$$
W^{1,p}_\Lambda(u^*T(T^*M))=\left\{
\begin{array}{ll}
\xi\in W^{1,p}(u^*T(T^*M)), \\
\xi(s,-1)\in T_{u(s,-1)}o_M,\,\xi(s,1)\in T_{u(s,1)}\nu^*N,\,s\in\mathbb{R},\\
\xi(s,0^-)\in T_{u(s,0^-)}\nu^*N,\,\xi(s,0^+)\in T_{u(s,0^+)}o_M,\,s\le0.
\end{array}
\right.
$$
We can see the operator $\overline{\partial}_{J,H}$ as a section of a Banach bundle
$$
\mathcal{E}\to\mathcal{P}^{1,p}(x_1^-,x_2^-;x^+),
$$
where a fiber over $u\in\mathcal{P}^{1,p}$ is
$$
\mathcal{E}_u=L^p(u^*T(T^*M)).
$$
Operator $\overline{\partial}_{J,H}$ is a Fredholm map. Linearization of this operator at its zero $u$ is given with
$$
E_u\xi:=D\overline{\partial}_{J,H}(u)\xi=\bigtriangledown_{\partial_su}\xi+J\bigtriangledown_{\partial_tu}\xi+\bigtriangledown_\xi J\frac{\partial u}{\partial t}-\bigtriangledown_\xi(JX_H(u)).$$
Fredholm property follows from the local elliptic estimate and the asymptotic properties (intersections $\nu^*N\cap\phi_{H_i}^1(o_M)$ are transversal). Similar as in Fredholm analysis for Lagrangian Floer homology we consider a set of almost complex structure on $T^*M$, $\mathcal{J}$, as a set of parameters and define a map
$$
\begin{aligned}
\mathcal{F}:\mathcal{P}&^{1,p}(x_1^-,x_2^-;x^+)\times\mathcal{J}\to\mathcal{E},\\
&(u,J)\mapsto\overline{\partial}_{J,H}(u).
\end{aligned}
$$
There exists a generic set $\mathcal{J}_{reg}\subset\mathcal{J}$ such that $E_u$ is onto for any $J\in\mathcal{J}_{reg}$ (explicit construction of Fredholm operator and analytic details will appear elsewhere). We conclude that for $J\in\mathcal{J}_{reg}$
$$
{\mathcal M}(x_1^-,x_2^-;x^+)=\overline{\partial}_{J,H}^{-1}(0)
$$
is a finite dimensional manifold with
$$
\dim{\mathcal M}(x_1^-,x_2^-;x^+)=\Ind\overline{\partial}_{J,H}.
$$
Compactness in the $C^\infty_{loc}$ topology of solutions of the Floer equation follows from Lemma 6.1, Proposition 6.2 in~\cite{AS} and from the fact that a non-constant $J$--holomorphic discs with boundary on $o_M\cup\nu^*N$ do not exist. Essentially, for compactness, there is no difference between strips with jumping boundary conditions considered by Abbondandolo and Schwarz and our jumping boundary condition on a slit of a pants. Removal of singularity (at a slit of a pants) of a $J$--holomorphic map with bounded energy follows from Proposition 6.5 in~\cite{AS}. Since $N$ is compact a sequence of pair--of--pants can break to something that is of the same kind plus holomorphic strip at the appropriate cylindrical end (see figure).
The boundary of 1--dimensional component of ${\mathcal M}(x,y;z)$ is a disjoint union
$$\begin{aligned}
\partial{\mathcal M}_{[1]}(x,y;z)=
&\bigcup_{x'\in CF_*(H_1)}{\mathcal M}(x,x';H_1)\times{\mathcal M}(x',y;z)\\
&\bigcup_{y'\in CF_*(H_2)}{\mathcal M}(y,y';H_2)\times{\mathcal M}(x,y';z)\\
&\bigcup_{z'\in CF_*(H_3)}{\mathcal M}(x,y;z')\times{\mathcal M}(z',z;H_3).
\end{aligned}$$
Gluing arguments in this situation are quite like various gluing arguments that appeared in different context (see~\cite{Sc1} or ~\cite{F1}). From pre-gluing of linearized Fredholm operator and an existence of an exact solution of Cauchy--Riemann equation near the glued strip and a pair--of--pants it follows that every element of a right--hand side is an element of a left--hand side.

\vskip5mm
\begin{center}\includegraphics[width=11cm,height=6cm]{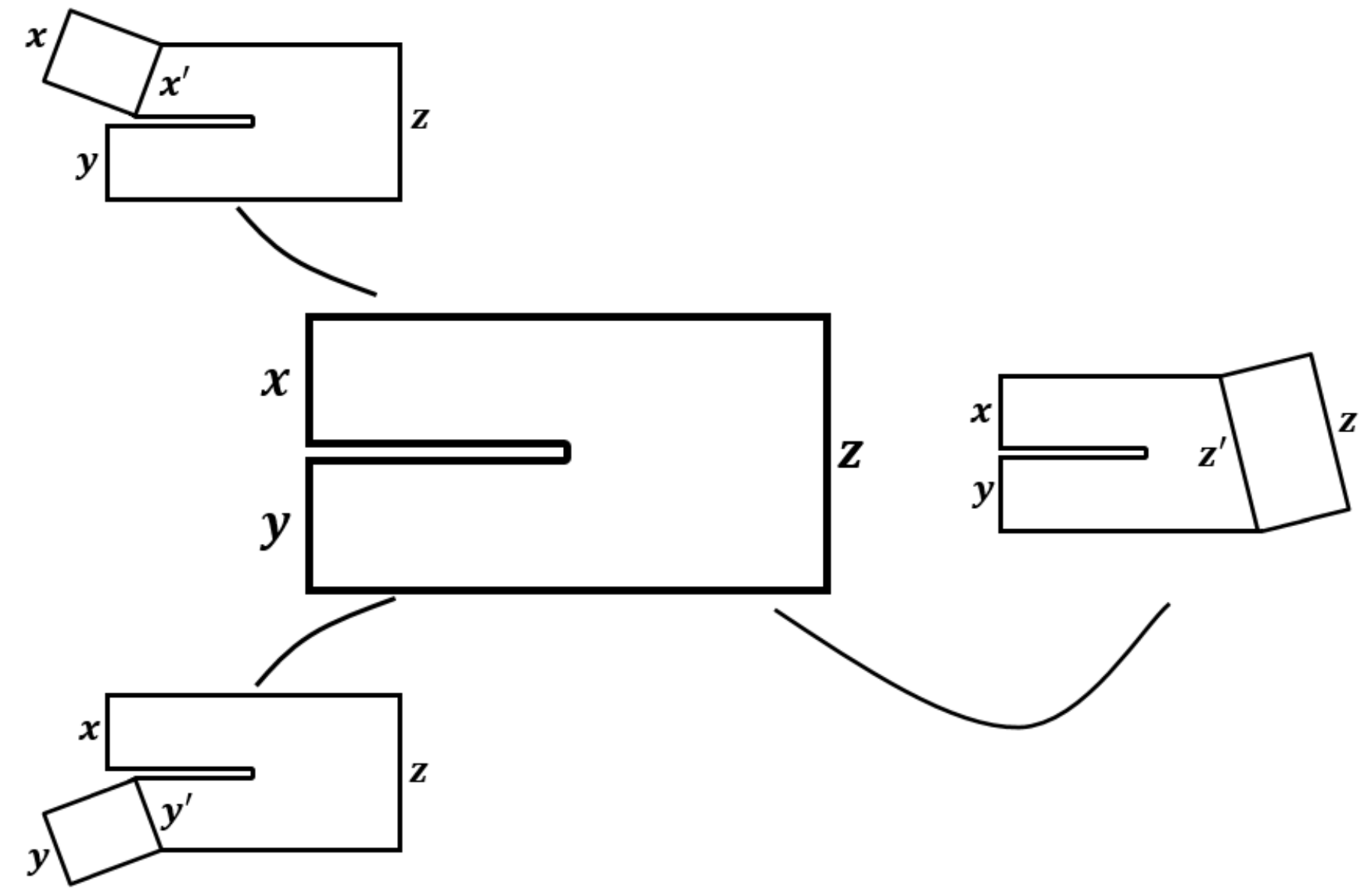}\\
Fig. 9. \textit{Boundary of ${\mathcal M}_{[1]}(x,y;z)$}
\end{center}
\vskip5mm
On generators of $CF_*$ we define
$$
x\ast y=\sum_z\sharp_2{\mathcal M}(x,y;z)\,z.
$$
Here, $\sharp_2{\mathcal M}(x,y;z)$ denotes the (modulo 2) number of elements of a zero--dimensional component of ${\mathcal M}(x,y;z)$. We extend the product $\ast$ by bilinearity on
$$CF_*(o_M,\nu^*N:H_1)\otimes CF_*(o_M,\nu^*N:H_2),$$
and conclude that $\ast$ commutes with the respective boundary operators and induces a product in homology.
\begin{rem}\label{alternative} We can define the product $\ast$ as a composition of a morphism
$$
m:HF_*(o_M,\nu^*N:H_1)\to HF_*(o_M,o_M:H_1')
$$
and the standard product in Lagrangian Floer homology
$$
\star:HF_*(o_M,o_M:H_1')\otimes HF_*(o_M,\nu^*N:H_2)\to HF_*(o_M,\nu^*N:H_3)
$$
(see~\cite{Aur} for the definition of the product or see proof of Proposition~\ref{bound} below). Morphism $m$ is defined on a chain level and counts holomorphic strips with jumping boundary conditions
$$
m(x)=\sum_{y}\sharp_2{\mathcal M}(x:H_1;y:H_1')\,y.
$$
Here ${\mathcal M}(x:H_1;y:H_1')$ denotes space of these holomorphic strips
$$\begin{aligned}&{\mathcal M}(x:H_1;y:H_1')=&\left\{u:{\mathbb R}\times[0,1]\to T^*M\left|\begin{array}{l}
 \partial_su+J(\partial_tu-X_{\overline{H}_1}\circ u)=0,\\
 u(s,0)\in o_M,\,s\in{\mathbb R},\\
 u(s,1)\in\nu^*N,\,s\le 0,\\
 u(s,1)\in o_M,\,s\ge0,\\
 u(-\infty,t)=x(t),\\
 u(+\infty,t)=y(t)
 \end{array}\right.\right\},\end{aligned}$$
and $\overline{H}_1$ is a homotopy that connects $H_1$ and $H_1'$. From standard gluing arguments it follows
$$\ast=\star\circ(m\otimes\mathbb{Id}).$$
The operation $m$ induces, via the PSS, the morphism
$$
H_*(N)\to H_*(M),
$$
while the operation $\star$ induces the morphism which is the action of $H_*(M)$ on $H_*(N)$.
\end{rem}



Thus, the product $\ast$ induces on $H_*(N)$, via the
PSS isomorphism, the operation given by composing the action of $H_*(M)$ on $H_*(N)$ and
the inclusion morphism $H_*(N)\to H_*(M)$.\\
\indent Now we can prove that conormal spectral invariants are subadditive with respect to $\ast$ product.\\

\noindent {\it Proof of Proposition~\ref{inv_ineq}:} Since a concatenation doesn't have to be a smooth function, we can find a Hamiltonian $H''$ that is regular, smooth and close enough to the concatenation $H\sharp H'$:
$$\|H''-H\sharp H'\|_{C^0}<\varepsilon.
$$
First step is to prove that the product $\ast$ defines a product on filtered complexes
$$
CF^{\lambda}_*(H)\times CF^{\mu}_*(H')\to CF_*^{\lambda+\mu+\varepsilon}(H''),
$$
for every $\varepsilon>0$ small enough.

Let us take smooth family of Hamiltonians $K:{\mathbb R}\times[-1,1]\times T^*M\to{\mathbb R}$ such that
$$K(s,t,\cdot)=\begin{cases}H(t+1,\cdot),&s\le -1,-1\le t\le0,\\
        H'(t,\cdot),&s\le -1,0\le t\le1,\\
        \frac{1}{2}H''(\frac{t+1}{2},\cdot),&s\ge1
        \end{cases}.$$
We can choose $K$ such that
$$
\bigg\|\frac{\partial K}{\partial s}\bigg\|\le\varepsilon,\;s\in[-1,1],
$$
and
$$
\frac{\partial K}{\partial s}=0,
$$
elsewhere.
Let us take $x\in CF^{\lambda}_*(H)$ and $y\in CF^{\mu}_*(H')$. Assume that there exists an element $u\in{\mathcal M}(x,y;z)$ for some $z\in CF_*(H'')$ ($u$ is a solution of an equation $\bar{\partial}_{K,J}(u)=0$). Then it holds
\begin{equation}\label{ff}
\begin{aligned}
0\le\int_{\Sigma}\bigg\|\frac{\partial u}{\partial s}\bigg\|^2\,ds\,dt&=\int_\Sigma\omega\bigg(\frac{\partial u}{\partial s},J\frac{\partial u}{\partial s}\bigg)\,ds\,dt\\
&=\int_\Sigma\omega\bigg(\frac{\partial u}{\partial s},\frac{\partial u}{\partial t}-X_K(u)\bigg)\,ds\,dt\\
&=\int_\Sigma u^*\omega-\int_\Sigma dK\bigg(\frac{\partial u}{\partial s}\bigg)\,ds\,dt.
\end{aligned}
\end{equation}
Using Stoke's formula we obtain
$$
\int_\Sigma u^*\omega=-\int x^*\lambda-\int y^*\lambda+\int z^*\lambda.
$$
Using the equality
$$
\int_\Sigma\frac{\partial}{\partial s}\big(K\circ u\big)\,ds\,dt=\int_\Sigma dK\bigg(\frac{\partial u}{\partial s}\bigg)\,ds\,dt+\int_\Sigma\frac{\partial K}{\partial s}(u)\,ds\,dt,
$$
and Stoke's formula again we get an estimate
$$
-\int_\Sigma dK\bigg(\frac{\partial u}{\partial s}\bigg)\,ds\,dt\le\int_0^1H(x(t),t)\,dt+\int_0^1H'(y(t),t)\,dt-\int_{0}^1H''(z(t),t)\,dt+4\varepsilon.
$$
Thus
$$
{\mathcal A}_{H''}(z)\le{\mathcal A}_H(x)+{\mathcal A}_{H'}(y)+4\varepsilon.
$$
From the the definition of operation $\cdot$ it easily follows
$$
l(\alpha\cdot\beta;o_M,\nu^*N:H'')\le l(\alpha;o_M,\nu^*N:H)+l(\beta;o_M,\nu^*N:H')+4\varepsilon.
$$
We know that spectral invariants are continuous with respect to the Hamiltonian (see \cite{O1}). If we pass to the limit as $\varepsilon\to0$ we get the triangle inequality
$$
l(\alpha\cdot\beta;o_M,\nu^*N:H\sharp H')\le l(\alpha;o_M,\nu^*N:H)+l(\beta;o_M,\nu^*N:H').
$$
\qed\\

Operation induced by $m$ provides inequality among spectral invariants. The proof of the following proposition is similar to the previous proof.
\begin{prop}
  For $\alpha\in H_*(N)\setminus\{0\}$ it holds
  $$l(\Phi(m(\Psi(\alpha)));o_M,o_M:H)\leq l(\alpha;o_M,\nu^*N:H).$$
\end{prop}

\noindent {\it Proof of Proposition~\ref{bound}:}
Let us take Morse functions $F:M\to{\mathbb R}$ and $f:N\to{\mathbb R}$. We want to define new type of a product on Morse homology:
$$
\bullet:HM_*(M:F)\otimes HM_*(N:f)\to HM_*(N:f).
$$
Let $p$ be a critical point of $F$ and $q,r$ critical points of $f$. We define $\bar{\mathcal {M}}(p,q;r)$ to be the set (see figure 10)
$$\begin{aligned}&\left\{(\Gamma,\gamma)\left|\begin{array}{l}
 \Gamma:(-\infty,0]\to M,\,\gamma:{\mathbb R}\to N,\\
 \dot{\Gamma}=-\nabla F(\Gamma),\,\dot{\gamma}=-\nabla f(\gamma),\\
 \Gamma(-\infty)=p,\,\gamma(-\infty)=q,\,\gamma(+\infty)=r,\\
 \Gamma(0)=\gamma(0)
 \end{array}\right.\right\}.\end{aligned}$$
\vskip5mm
\begin{center}\includegraphics[width=8cm,height=2.5cm]{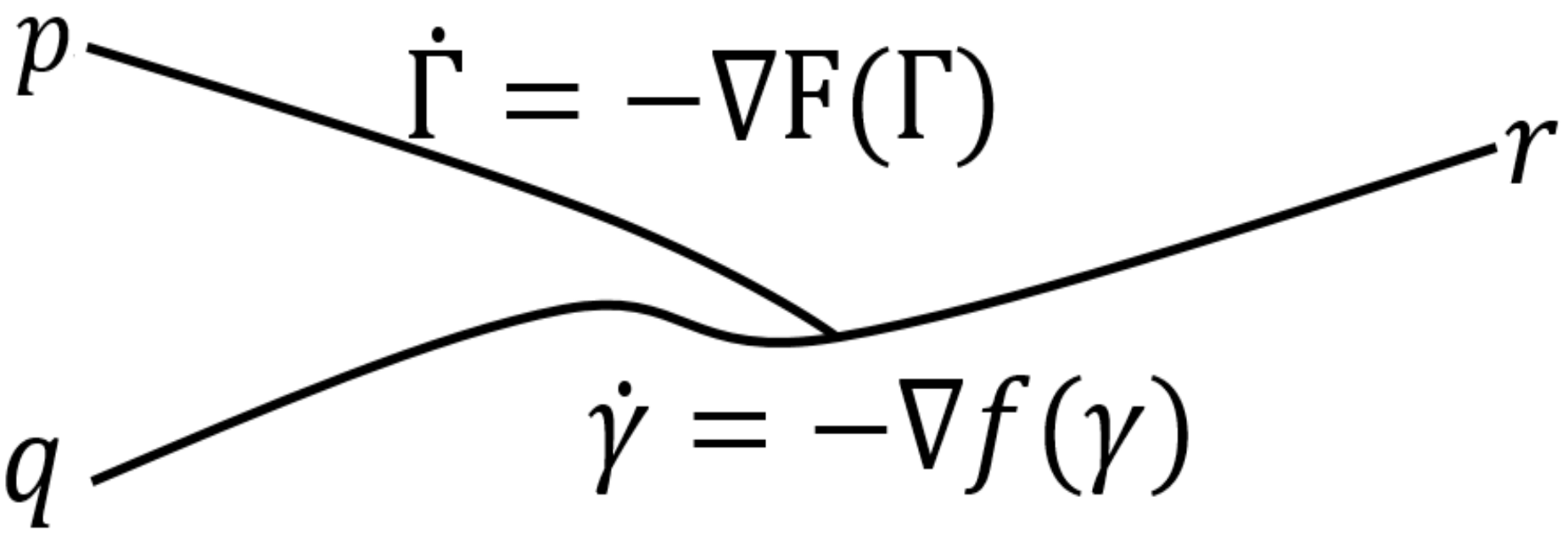}\\
Fig. 10. \textit{$\bar{\mathcal {M}}(p,q;r)$}
\end{center}
\vskip5mm
Let us compute the dimension of $\bar{\mathcal {M}}(p,q;r)$. We define a map
$$
\begin{aligned}
ev:W&^{u}(p,F)\times{\mathcal M}(q,r)\to M\times N,\\
&(\Gamma,\gamma)\mapsto(\Gamma(0),\gamma(0)).
\end{aligned}
$$
For generic choices, $ev$ is transversal to a submanifold
$$
\Delta_N=\{(x,x)\,|\,x\in N\}\subset M\times N
$$
and
$$
\bar{\mathcal {M}}(p,q;r)=ev^{-1}(\Delta_N).
$$
Simple computation gives
$$
\dim\bar{\mathcal {M}}(p,q;r)=m_F(p)+m_f(q)-m_f(r)-\dim M.
$$
If we denote by $\bar{n}(p,q;r)$ the number of elements of a zero--dimensional component of $\bar{\mathcal M}(p,q;r)$ we can define a product $\bullet$ $$
p\bullet q=\sum_{r}\bar{n}(p,q;r)\,r.
$$
This map $\bullet$ agrees with the boundary operator and it induces a product in homology.

Specially, we can take a Morse function $F$ that has a unique critical point $p$ of an index $m_F(p)=\dim M$ (unique maximum). A Morse homology class of this point represents the fundamental class in $H_{\dim M}(M)$. Then $\bar{n}(p,q;r)$ counts number of pairs $(\Gamma,\gamma)$, where $\gamma$ is trajectory that connects critical points $q$ and $r$ such that $m_f(q)=m_f(r)$. Number of such trajectories $\gamma$ is 0 if $q\neq r$ and 1 if $q=r$ (constant trajectory). Now we want to find the number of negative gradient trajectories $\Gamma$ that leaves global maximum $p$ and hits a point $q=r$. We can pick a generic function $f$ such that its critical points belong to $W^{u}(p)$. Since $q=r\in W^{u}(p)$ such trajectory $\Gamma$ exists and is unique. We conclude that the multiplication with a class $[p]$ induces the identity on the homology:
$$
{\mathbb {Id}}=[p]\bullet: HM_*(N:f)\to HM_*(N:f).
$$

Finally, in order to prove the boundness of spectral invariants, we need to describe a construction of the standard product in Lagrangian Floer homology:
 $$ \star: HF_*(o_M,\nu^*N:H_1)\otimes HF_*(o_M,o_M:H_2)\to HF_*(o_M,\nu^*N:H_3).
 $$
Note that $HF_*(o_M,o_M:H_2)$ is Floer homology for conormal bundle in a special case when $M=N$.
Similarly to the construction of the product $\ast$, we consider the space of perturbed holomorphic maps on a Riemannian surface $\Sigma$ but with different boundary conditions. For $x\in CF_*(o_M,o_M:H_2)$, $y\in CF_*(o_M,\nu^*N:H_1)$ and $z\in CF_*(o_M,\nu^*N:H_3)$ we define
$$\begin{aligned}&\bar{\mathcal M}(x,y;z)=&\left\{u:\Sigma\to T^*M\left|\begin{array}{l}
 \partial_su_j+J(\partial_tu_j-X_{\rho_jH_j}\circ u_j)=0,\,j=1,2,3,\\
 \partial_su+J\partial_tu=0,\,\mbox{on}\;\Sigma_0:=\Sigma\setminus(\Sigma_1\cup \Sigma_2\cup \Sigma_3),\\
 u(s,-1)\in o_M,\,u(s,1)\in\nu^*N,\,s\in{\mathbb R},\\
 u(s,0^-),u(s,0^+)\in o_M,\,s\le 0,\\
 u_1(-\infty,t)=x(t),\\
 u_2(-\infty,t)=y(t),\\
 u_3(+\infty,t)=z(t)
 \end{array}\right.\right\}.\end{aligned}$$
With
$$
x\star y=\sum_{z}\bar{n}(x,y;z)\,z
$$
we define a map on a chain complex that defines a product in homology. Here $\bar{n}(x,y;z)$ denotes the (modulo 2) number of elements of a zero--dimensional component of $\bar{\mathcal M}(x,y;z)$. (Similar type of product is defined in~\cite{DKM}. They use it to compare spectral invariants in Lagrangian and Hamiltonian Floer theory.) Using the standard cobordism arguments it follows
$$
\Psi^{\nu}(\alpha\bullet\beta)=\Psi^{\nu}(\alpha)\star\Psi^{o}(\beta),
$$
for $\alpha\in HM_*(N:f_1)$ and $\beta\in HM_*(M:f_2)$. Here $\Psi^{\nu}$ and $\Psi^{o}$ denote PSS isomorphisms
$$
\begin{aligned}
&\Psi^{\nu}:HM_*(N:f)\to HF_*(o_M,\nu^*N:H),\\
&\Psi^{o}:HM_*(M:f)\to HF_*(o_M,o_M:H).
\end{aligned}
$$
Using the same argument as in the proof of Proposition ~\ref{inv_ineq} one can prove that it holds
$$
l(\alpha\bullet\beta;o_M,\nu^*N:H\sharp H')\le l(\alpha;o_M,\nu^*N:H)+l(\beta;o_M,o_M:H'),
$$
for all $\alpha\in H_*(N)$ and $\beta\in H_*(M)$ such that $\alpha\bullet\beta\neq0$.
Specially, if we take $\beta=[M]$ and $H=0$ we obtain an inequality
$$
l(\alpha;o_M,\nu^*N:0\sharp H')\le l(\alpha;o_M,\nu^*N:0)+l([M],o_M,o_M:H'),
$$
that holds for all $\alpha\in H_*(N)\setminus\{0\}$. Since spectral invariants are continuous and they belong to the spectrum of Hamiltonian $H$ it follows that
$$
l(\alpha;o_M,\nu^*N:0)=0.
$$
The concatenation $0\sharp H'$ is just a reparametrization of $H'$ and it doesn't change Hamiltonian orbits, Floer strip and spectral invariants. Thus
$$l(\alpha;o_M,\nu^*N:0\sharp H')=l(\alpha;o_M,\nu^*N:H').$$
We conclude that spectral invariants of non--zero homology classes are bounded from above:
$$l(\alpha;o_M,\nu^*N:H')\le l([M],o_M,o_M:H').
$$
\qed
\smallskip


\begin{thebibliography}{11}

\bibitem{AS} A. Abbondandolo, M. Schwarz, {\it Floer homology of cotangent bundles and the loop
product},  Geom. Topol. 14(3) (2010), 1569-–1722.

\bibitem{Alb} P. Albers, {\it A Lagrangian Piunikhin-Salamon-Schwarz morphism and two comparison
homomorphisms in Floer homology }, Int. Math. Res. Not. IMRN 2008, no. 4, 56pp.

\bibitem{Aur} D. Auroux, {\it A Beginner’s Introduction to Fukaya Categories}, arXiv:1301.7056 (2013).

\bibitem{BC} P. Biran, O. Cornea, {\it Quantum structures for Lagrangian submanifolds},
http://arxiv.org/pdf/0708.4221.

\bibitem{CTL} K. Cieliebak, T. Ekholm, J. Latschev,{\it Compactness for holomorphic curves with
switching Lagrangian boundary conditions}, J. Symplectic Geom. 8 no. 3 (2010), 267--298.

\bibitem{DKM} J. \Dj ureti\'c, J. Kati\'{c}, D. Milinkovi\'{c}, {\it Comparison of spectral invariants
in Lagrangian and Hamiltonian Floer theory}, Filomat 30:5 (2016), 1161–1174.

\bibitem{F1} A. Floer, {\it Morse theory for Lagrangian intersections}, J. Differential Geom. 28 (1988), 513--547.

\bibitem{F2} A. Floer, {\it Symplectic fixed points and holomorphic spheres}, Comm. Math. Phys., 120 (1989), 575--611.

\bibitem{Urs} U. Frauenfelder, {\it Floer homology of symplectic quotients and the Arnold--Givental conjecture},
PhD thesis, ETH Z\"{u}rich, 2003.

\bibitem{Urs1} U. Frauenfelder, {\it Gromov convergence of pseudoholomorphic discs},
Journal of Fixed Point Theory and Application, Volume 3 (2008), Number 2, 215--271.

\bibitem{FOOO} K. Fukaya, Y.-G. Oh, H. Ohta, K. Ono, {\it Lagrangian
intersection Floer theory}, Kyoto University preprint.

\bibitem{G} M. Gromov, {\it Pseudo holomorphic curves in symplectic manifolds}, Invent. Math. 82 (1985),
307–-347.

\bibitem{HLS} V. Humili\`{e}re, R. Leclercq, S. Seyfaddini, {\it Reduction of symplectic homeomorphisms}, arXiv:1407.6330v2 (2014).

\bibitem{K} J. Kati\'{c}, {\it Compactification in mixed moduli spaces in Morse--Floer theory},
Rocky Mountain Journal of Mathematics, 38 (2008), 923--939.

\bibitem{KM} J. Kati\'{c}, D. Milinkovi\'{c}, {\it Piunikhin--Salamon--Schwarz isomorphism for Lagrangian intersections}, Diff. Geom. and its Appl., 22 (2005), 215--227.

\bibitem{KMS} J. Kati\'{c}, D. Milinkovi\'{c}, T. Sim\v{c}evi\'{c}, {\it Isomorphism between Morse and Lagrangian Floer cohomology rings}, Rocky Mountain Journal of Mathematics, 41, No. 3 (2011), 789--811.

\bibitem{MS} D. McDuff, D. Salamon, {\it J-holomorphic Curves and
Quantum Cohomology}, AMS, University Lecture Series 6, 1994.

\bibitem{MS1} D. McDuff, D. Salamon, {\it J-holomorphic Curves and Symplectic Topology},
American Mathematical Society Colloquium Publications, vol. 52, AMS, Providence, RI, 2004.

\bibitem{M1} D. Milinkovi\'c, {\it Morse homology for generating functions of Lagrangian submanifolds},
Trans. Amer. Math. Soc., Vol. 351, No. 10 (1999), 3953--3974.

\bibitem{M2} D. Milinkovi\'c, {\it On equivalence of two constructions of invariants of Lagrangian submanifolds},
Pacific J. Math., Vol. 195, No. 2 (2000), 371--475.

\bibitem{M} J. Milnor, {\it Lectures on the h-cobordism Theorem},
Princeton University Press, 1963.

\bibitem{MNZ} A. Monzner, N. Vichery, F. Zapolsky, {\it Partial quasi-morphisms and quasi-states on cotangent
bundles, and symplectic homogenization}, Journal of Modern Dynamics, Issue 2 (2012), 205–-249.

\bibitem{O1} Y.-G. Oh, {\it Symplectic topology as the geometry of action
functional I -- relative Floer theory on the cotangent bundle}, J.
Differential Geom. 45 (1997), 499--577.

\bibitem{O2} Y.-G. Oh, {\it Symplectic topology as the geometry of action
functional, II - pants product and cohomological invariants}, Comm. Anal. Geom. 7 (1999), 1--55.

\bibitem{PSS} S. Piunikhin, D. Salamon, M. Schwarz, {\it Symplectic
Floer--Donaldson theory and quantum cohomology}, in: Contact and
symplectic geometry, Publ. Newton Instit. 8, Cambridge Univ. Press,
Cambridge (1996), pp. 171--200.

\bibitem{P} M. Po\'{z}niak, {\it Floer homology, Novikov rings and clean
intersections}, Ph.D. thesis, University of Warwick, 1994.

\bibitem{RS1} J. Robbin, D. Salamon, {\it The Maslov index for paths},
Topology 32 (1993), 827--844.

\bibitem{RS2} J. Robbin, D. Salamon, {\it The spectral flow and the
Maslov index}, Bull. London Math. Soc. 27 (1995), 1--33.

\bibitem{Sa1} D. Salamon, {\it Morse theory, the Conley index and Floer
homology}, Bull. Lond. Math. Soc. 32 (1990), 113--140.

\bibitem{Sa} D. Salamon, {\it Lectures on Floer homology}, IAS Park
City Math. Series, AMS Vol 7, 1999.

\bibitem{FelS} F. Schm\"{a}schke, {\it Floer homology of Lagrangians in clean intersection},
arXiv:1606.05327, 2016.

\bibitem{Sc1} M. Schwarz, {\it Morse Homology}, Birkh\"{a}user, 1993.

\bibitem{Sim} T. Sim\v{c}evi\'c, {\it A Hardy Space Approach to Lagrangian Floer gluing},
PhD thesis, ETH Z\"{u}rich, 2014.

\bibitem{V} C. Viterbo, {\it Symplectic topology as the geometry of
generating functions}, Math. Ann. 292 no. 4 (1992), 685--710.


\end{thebibliography}
\end{document}